\documentclass[a4paper,10pt,appendixprefix=true]{scrbook}
\usepackage[T1]{fontenc}
\usepackage[utf8]{inputenc}
\usepackage[mono=false]{libertine}
\usepackage[most]{tcolorbox}
\usepackage[margin=2cm]{caption}
\usepackage[shortlabels]{enumitem}
\usepackage[style=numeric,sorting=ymdnt,date=year]{biblatex}
\usepackage{amsthm,amsmath,amssymb,stmaryrd,mathtools,bm}
\usepackage{caption,booktabs,floatrow,graphicx,xcolor,svg,pdfpages}
\usepackage[pdfusetitle]{hyperref}

\date{}
\title{An Exploration of Graph Pebbling}
\author{Herman Bergwerf}
\hypersetup{
  hidelinks,
  pdfauthor={Herman Bergwerf},
  pdfstartview=Fit,
  pdfpagemode=UseOutlines,
  pdfpagelayout=TwoPageRight,
  bookmarksnumbered=true
}

\DeclareSortingTemplate{ymdnt}{
  \sort{\field{presort}}
  \sort{\field{sortyear}\field{year}}
  \sort{\field{month}}
  \sort{\field{day}}
  \sort{\field{sortname}\field{author}}
  \sort{\field{sorttitle}\field{title}}
}

\renewbibmacro*{name:andothers}{
  \ifboolexpr{test {\ifnumequal{\value{listcount}}{\value{liststop}}} and test \ifmorenames}{%
    \ifnumgreater{\value{liststop}}{1}{\finalandcomma}{}%
    \andothersdelim\bibstring[\emph]{andothers}}{}}

\DeclareMathSizes{10}{10}{6}{4}

\counterwithout{figure}{chapter}
\counterwithout{table}{chapter}
\DeclareCaptionLabelFormat{custom}{\textbf{#1 #2}}
\DeclareCaptionFormat{custom}{#1#2#3}
\DeclareCaptionLabelSeparator{custom}{\:--\:}
\newcommand{\SimpleFigure}[2]{%
\begin{figure}[h]
\vspace{2mm}%
\includesvg{#1}%
\vspace{2mm}%
\caption{#2}
\end{figure}%
}

\newenvironment{Figure}[1]{%
\figure[h]
\vspace{2mm}%
\includesvg{#1}%
\vspace{2mm}%
}{\endfigure}

\newtheorem{question}{Question}
\newtheorem{conjecture}{Conjecture}
\newtheorem{definition}{Definition}[chapter]

\newtheorem{lemma}[definition]{Lemma}
\newtheorem{theorem}[definition]{Theorem}
\newtheorem{corollary}[definition]{Corollary}
\newtheorem*{theorem*}{Theorem}

\newcommand{\Erdos}{Erdős}
\newcommand{\ie}{\textit{i.e.}~}

\newcommand{\witregel}{\vspace{5mm}}

\newcommand{\fF}{\mathcal{F}}
\newcommand{\Nat}{\mathbb{N}}
\newcommand{\Frac}{\mathbb{Q}^+}
\newcommand{\dom}{\mathrm{dom}}
\newcommand{\defeq}{\vcentcolon=}
\renewcommand{\implies}{\Rightarrow}

\newcommand{\If}{\mathrm{if}}
\newcommand{\Else}{\mathrm{otherwise}}
\newcommand{\Forall}[2]{\forall #1:#2}
\newcommand{\Exists}[2]{\exists #1:#2}
\newcommand{\SmallCount}[1]{\raisebox{0.8pt}{\scalebox{0.6}{$\#$}}{#1}}
\newcommand{\Count}[1]{\raisebox{1.4pt}{\scalebox{0.8}{$\#$}}{#1}}
\newcommand{\Set}[2]{\left\{#1~\middle\vert~#2\right\}}

\newcommand{\fplus}{+}
\newcommand{\fmin}{-}
\newcommand{\fprod}{\cdot}
\newcommand{\supp}[1]{\mathrm{supp}(#1)}
\newcommand{\size}[1]{\left|#1\right|}
\newcommand{\rsize}[2]{\mathsf{r}_{#1}\mathopen{}\left(#2\right)\mathclose{}}
\newcommand{\suppsize}[1]{\mathsf{s}\raisebox{2pt}{\scalebox{0.9}{\normalfont\#}}(#1)}
\newcommand{\bincount}[1]{\mathrm{count}(#1)}

\newcommand{\edge}[3]{({#2},{#3})\in E({#1})}
\newcommand{\weight}{\omega}

\newcommand{\diam}{\mathrm{diam}}
\newcommand{\gprodbox}{\raisebox{1pt}{\scalebox{0.6}{$\square$}}}
\newcommand{\gprod}{\mathbin{\hspace{2pt}\gprodbox\hspace{2pt}}}
\newcommand{\completegraph}[2]{K_{#1}^{(#2)}}
\newcommand{\pathgraph}[2]{P_{#1}^{(#2)}}

\newcommand{\instargraph}[2]{R_{#1}^{(#2)}}
\newcommand{\arrowgraph}[1]{{\downarrow^{(#1)}}}
\newcommand{\arrowfrom}{0}
\newcommand{\arrowto}{1}
\newcommand{\Hom}[3]{#3:#1\to #2}

\newcommand{\rtcto}{\to^\ast}

\newcommand{\solv}[3]{\mathrm{Solv}_{#1,#2}(#3)}
\newcommand{\inflow}[2]{\mathrm{in}_{#1}(#2)}
\newcommand{\outflow}[2]{\mathrm{out}_{#1}(#2)}
\newcommand{\outfloww}[2]{\mathrm{out}^\star_{#1}(#2)}
\newcommand{\excess}[2]{x_{#1}(#2)}

\newcommand{\pebb}[1]{\pi\mathopen{}\left(#1\right)\mathclose{}}
\newcommand{\npebb}[2]{\pi_{#1}\mathopen{}\left(#2\right)\mathclose{}}
\newcommand{\optpebb}[1]{\pi_\mathrm{opt}\mathopen{}\left(#1\right)\mathclose{}}
\newcommand{\taupebb}[3]{\tau_{#1,#2}\mathopen{}\left(#3\right)\mathclose{}}

\begin{document}
\pagenumbering{gobble}
\includepdf{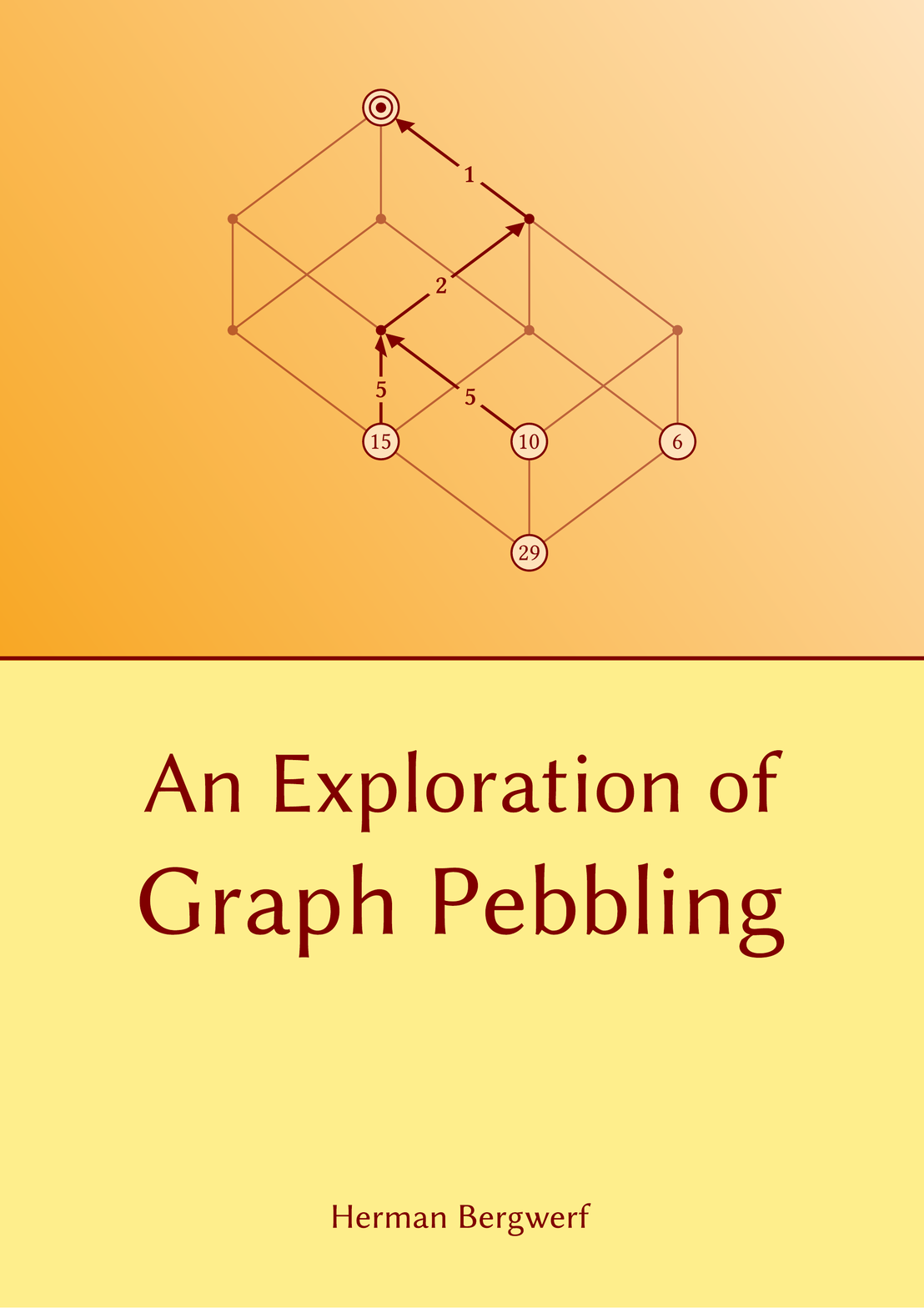}
\cleardoublepage
\pagenumbering{arabic}

\chapter*{Preface}
I wrote this thesis for my master's degree at Radboud University. It is the culmination of a five-year journey into mathematics and formal logic, which was initiated by a side project where I attempted to build my own theorem prover. Here I want to briefly tell a story about this journey.

\begin{center}
\includesvg[width=5mm]{fig/misc/asterism.svg}
\end{center}

When I first went to university eight years ago I had little interest in mathematics. I wanted to understand more about the machinery of life and decided to join the \emph{Nanobiology} program in Delft, which offered a glimpse into the complex world of low-level biological systems. During this time I worked on a side project called \emph{QEDb}, which aimed to develop a structured model and an attractive interface for writing and storing mathematical derivations. Initially I had solutions to physics problems in mind, such as the exercises we had to do for our classes. The derivation database would be a way for students to share their solutions in a rigorous way, and to learn from the solutions of others.
I worked on the QEDb project for several years, but, as I learned more about mathematical proofs, I kept finding things that required fundamental changes to my program. At first it only reasoned with equality and rewriting, then I implemented deduction rules, and at some point I learned about quantifiers and realized they had to be expressible in the system too. Eventually I abandoned the project, but it had sparked a definite interest in mathematics and formal systems.

\witregel

I had known about the Coq theorem prover for some time, but I didn't have a clue about how it worked or how I could learn it. I thought that the most effective way to achieve my goal was to build my own formal system. A course from Robbert Krebbers, where I first learned how to use Coq, changed this. I realized that Coq solved many problems I had encountered, and that there was a great deal for me to learn about theorem provers. This was an important motivation to study theoretical computer science in Nijmegen, where courses were offered about type systems, proof assistants, formal semantics, and other interesting topics.
I was also increasingly fascinated by purely mathematical problems. Some theorems, like Gödel's incompleteness theorem, or the decidability of linear integer arithmetic, were surrounded by an air of mystery, and I \emph{had} to know more about these things; I wanted to know how such results were even \emph{possible}. I particularly enjoyed the courses of Wim Veldman, whose classes were full of interesting stories and cliffhangers. His course about the independence of the continuum hypothesis from ZFC set theory is perhaps the most sophisticated and intriguing bit of theory I have ever studied.

\pagebreak
For my thesis I wanted to dive into a similarly mystifying topic, and do something that is not directly related to computer programming. When I asked Wim Veldman how to find an interesting problem to work on, he told me to read about a subject I liked and just ask myself questions about it, no matter how simple. He recommended the book \emph{The Strange Logic of Random Graphs} by Joel Spencer, which among other things treats Ehrenfeucht-Fraïssé games. In a related article I encountered the term \emph{pebbling game} that, quite by accident, brought me to the topic of this thesis (which is actually completely unrelated to what that article was talking about). For over half a year I puzzled feverishly over various problems related to graph pebbling, eager to find a gold nugget that would reveal the path to a new solution. In the end I wasn't able to solve any of the open problems, but I got the opportunity to develop my own intuition, and address some small questions along the way. In this thesis, I want to show what a rich and fascinating topic graph pebbling is. I hope it will inspire others to study it too!

\begin{center}
\vspace{5mm}
\includesvg[width=5mm]{fig/misc/asterism.svg}
\vspace{5mm}
\end{center}

There is a number people that played a role in this process and who I am especially grateful to. First of all I want to thank Wieb Bosma, who despite being unfamiliar with graph pebbling was prepared to supervise my thesis. We spent a lot of time discussing various results, and Wieb was always open to exploring new directions. He gave me a lot of detailed feedback that sharpened the style and quality of my writing. Without him I would not have reached this point. I am also very thankful to Wim Veldman, who through his classes inspired a certain way of thinking and writing about mathematics. His enthusiasm and rich use of metaphors were very motivating, and I never so enjoyed sitting through lectures. His course about intuitionistic mathematics revealed a whole new way of looking at logic. Finally, I also want to thank Robbert Krebbers, from whom I learned developing efficient Coq proofs, and who has, over the past years, answered many of my emails about Coq.

\witregel

\begin{flushright}
\itshape
Herman Bergwerf\\
April 2023, Nijmegen
\end{flushright}

\protect\setcounter{tocdepth}{0}
{\Large\tableofcontents}

\chapter*{Notation}
\begin{itemize}[--,labelsep=3mm,leftmargin=*]
\item The set of \emph{natural numbers}, including the number zero, is denoted by $\Nat$.
\item The set of \emph{non-negative rational numbers} is denoted by $\Frac$.
\item The \emph{number of elements} of a set $S$ is written as $\Count{S}$.
\item The \emph{domain} of a function $f$ is written as $\dom(f)$.
\item The \emph{support} of a real-valued function $f$ is defined as~\,$\supp{f}\defeq\Set{x\in\dom(f)}{f(x)\neq 0}$.
\item The \emph{size} of a real-valued function $f$ is defined as $\size{f}\defeq\sum_{x\in\dom(f)}f(x)$.
\item The \emph{element-wise addition} of two real-valued functions $f$ and $g$ is written as $f\fplus g$.
\item The \emph{element-wise subtraction} of two real-valued functions $f$ and $g$ is written as $f\fmin g$.
\item The \emph{element-wise multiplication} of two real-valued functions $f$ and $g$ is written as $f\fprod g$.
\item The \emph{complete graph} on $n$ vertices is denoted by $K_n$.
\item The \emph{cycle graph} on $n$ vertices is denoted by $C_n$.
\item The \emph{path graph} on $n$ vertices is denoted by $P_n$.
\end{itemize}

\chapter{Introduction}
\label{chapter:introduction}

The topic of this treatise is a combinatorial technique called \emph{Graph Pebbling}. Before we start investigating this, it is worth considering where such a thing comes from. In the early 20th century one of the great unsolved problems in mathematics was the \emph{Entscheidungsproblem}, which asked if there exists a procedure to determine if a given sentence in first-order logic is universally true. One of the people who worked on this was Frank Ramsey\footnote{For a brief outline of Ramsey's life, see \cite{Mellor1983}.}, who found such a procedure for a restricted class of logical formulae\footnote{``On a Problem of Formal Logic'' (1930)}. It is now well known that a general decision procedure for first-order logic does \emph{not exist}, which was shown by both Church\footnote{``An Unsolvable Problem of Elementary Number Theory'' (1936)} and Turing\footnote{``On Computable Numbers, with an Application to the Entscheidungsproblem'' (1936)} using formal systems to express \emph{algorithms}. Ramsey didn't witness this result; he passed away at just 26. In his work on the Entscheidungsproblem he introduced a curious theorem about graphs that carried an intriguing concept; certain orderly substructures are present in any sufficiently large structure, \emph{no matter how disordered}. This result inspired a whole mathematical field that is now known as \emph{Ramsey Theory}.

\witregel

The renowned mathematician Paul \Erdos~helped to popularize Ramsey Theory, and worked on many problems in this new field. In 1961 he showed\footnote{``Theorem in the Additive Number Theory'' (1961)}, together with Ginzburg and Ziv, that every integer sequence of length $2n-1$ has a subsequence of $n$ integers that together sum to a multiple of $n$, \ie~zero modulo $n$. This seminal work inspired research into other \emph{zero-sum} problems. Communication between \Erdos~and the mathematician Paul Lemke resulted in the theorem stated below, which we shall call the \emph{\Erdos-Lemke conjecture}.

\begin{theorem}
\label{thm:erdos_lemke_conj}
\emph{(\Erdos-Lemke conjecture)} Given positive integers $n$, $d$, and $a_1,a_2,\dots,a_d$ with $d\,|\,n$ and $a_i\,|\,n$ for $i=1,2,\dots,d$, there is a non-empty subset $S$ of $\{1,2,\dots,d\}$ such that:
$$\sum_{i\in S}a_i\equiv 0~(\mathrm{mod}~d)\qquad\text{and}\qquad\sum_{i\in S}a_i\le n.$$
\end{theorem}

A proof of this conjecture was published by Lemke and Kleitman in \cite{Lemke1989} a few years after its conception. Their proof, which uses induction on the number of prime factors of $d$, is quite complicated and difficult to divide into smaller steps. Fan Chung developed an alternative proof in \cite{Chung1989} based on a game of moving pebbles between the vertices of a graph, an idea originating from Lagarias and Saks (who are absent in the literature). In this proof the various reasoning steps are more clearly separated, and most work is spent analyzing the pebbling game on a hypercube. Chung raised some interesting questions about the pebbling game on other graphs, inspiring a string of further research into what is now called Graph Pebbling\footnote{For a recent survey of graph pebbling research, see \cite{Hurlbert2021}.}.

\begin{figure}[h]
\includegraphics[width=.5\textwidth]{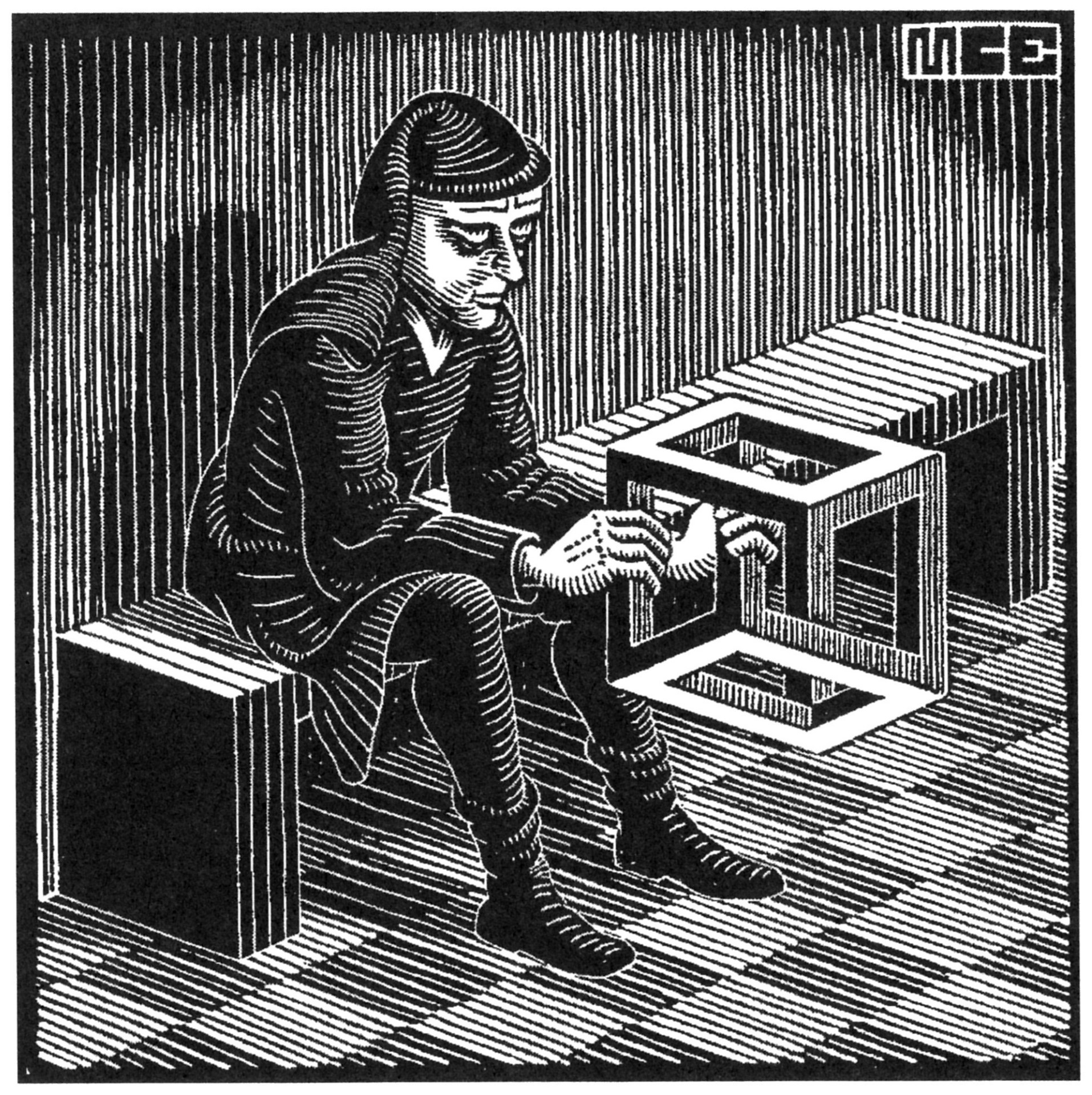}
\caption{\emph{Man with Cuboid} by M.C. Escher (1958)}
\end{figure}

Let's consider an instance of this pebbling game. Suppose I distribute seven pebbles over the vertices of a cube, where it is possible to put multiple pebbles on a single vertex. Starting from this initial configuration of pebbles, you are allowed to make changes with the following rule: You may add \emph{one} pebble to any vertex in exchange for removing \emph{two} pebbles from an adjacent vertex, resulting in a net loss of one pebble. This is called a \emph{pebbling step}. Your objective is to put one pebble on a designated empty vertex using pebbling steps. Can I pick a certain initial configuration with an empty vertex for which you are unable to achieve the objective? And what happens if I distribute eight pebbles instead of seven?


\witregel

We are going to explore this game in detail in the coming chapters. In the next chapter we will first define some graph pebbling concepts formally, and get more familiar with them by looking at pebbling in trees and diameter-2 graphs. We will also discuss some questions that reveal the complexity of our topic. In Chapter~\ref{chapter:weight} we explore a technique that uses linear inequalities to determine if an initial pebble configuration can be solved by a specific strategy, yielding a number of elegant proofs and a method to apply linear programming. In Chapter~\ref{chapter:flow} we think about the direction in which pebbles need to be moved, and we prove that back-and-forth movement between vertices is never a necessity. In Chapter~\ref{chapter:hypercubes} we review Chung's results about pebbling on a hypercube, and we describe a more detailed proof with a more general induction argument. Finally, in Chapter~\ref{chapter:zero-sums}, we will prove the \Erdos-Lemke conjecture via graph pebbling. Our investigation of this result was driven by a desire to describe a formal proof using the \emph{Coq Proof Assistant}. Details of the completed formalization are described in the appendix.


%
%
%

\chapter{Pebbling Numbers}
\label{chapter:pebbling}

Contrary to most of the literature about graph pebbling, we are going to describe the pebbling game using directed graphs and edge weights. This will be particularly useful for our proof of the \Erdos-Lemke conjecture.

\begin{definition}
A simple graph $G$ is a pair $(V,E)$ of a set of vertices $V$ and a set of edges $E$, where each edge $\{u,v\}\in E$ is a set of two distinct vertices. The vertices of $G$ are written as $V(G)$ and the edges as $E(G)$.
\end{definition}

\begin{definition}
An edge-weighted digraph $G$ is a triplet $(V,E,\weight)$ of a set of vertices $V$, a set of edges $E\subseteq V\times V$, and an edge-weight function $\weight:E\to\Nat$ that assigns a weight to every edge. The vertices of $G$ are written as $V(G)$, the edges as $E(G)$, and the edge-weight function as $\weight_G$.\linebreak
An edge $(u,v)\in E$ begins at $u$ and ends at $v$.
\end{definition}

\noindent We refer to edge-weighted digraphs simply as \emph{graphs}. Sometimes graphs with undirected edges are given, in which case you may assume that all edges are bidirectional. We define a special notation to convert a simple graph, such as the complete graph $K_n$, to an edge-weighted digraph where all edges are bidirectional, and the weight of every edge is $k$:

\begin{definition}
Let $G$ be a simple graph and let $k\in\Nat$.
The graph $G^{(k)}$ is defined as $(V,E,\weight)$ where $V=V(G)$, $E=\Set{(u,v)}{\{u,v\}\in E(G)}$, and $\weight(u,v)=k$ for all $u,v\in V(G)$.
\end{definition}

Throughout this text, we will encounter a special way to construct larger graphs from smaller graphs known as a \emph{Cartesian graph product}. The Cartesian product of two graphs $G$ and $H$ is written as $G\gprod H$. The box symbol is a reference to the orthogonality of the construction; the product of two paths is a grid, and the product of two cycles is a torus.

\begin{definition}
Let $G$ and $H$ be graphs.
The product $G\gprod H$ is defined as $(V,E,\weight)$ where:
\begin{align*}
V&= V(G)\times V(H)\\[2mm]
E&=\Set{\begin{aligned}((u_1,u_2),(v_1,v_2))\\\in V\times V\end{aligned}}{\begin{aligned}
&(\edge{G}{u_1}{v_1}\land u_2=v_2)~\lor\\
&(\edge{H}{u_2}{v_2}\land u_1=v_1)
\end{aligned}}\\[2mm]
\weight((u_1,u_2),(v_1,v_2))&=\begin{cases}
\weight_G(u_1,v_1)&\If~\edge{G}{u_1}{v_1}\\
\weight_H(u_2,v_2)&\If~\edge{H}{u_2}{v_2}
\end{cases}
\end{align*}
\end{definition}

\pagebreak
\SimpleFigure{fig/example/Cartesian_product_torus.svg}{The product of two cycles is a torus.}

Graph pebbling starts with a distribution of pebbles over the vertices of a graph, called the \emph{initial configuration}. This configuration can be altered with \emph{pebbling steps}: If an edge $(u,v)$ has weight $k$, and $u$ has at least $k$ pebbles, then we may remove $k$ pebbles from $u$ and add one pebble to~$v$. Essentially all of the literature about graph pebbling defines pebbling steps to require a standard amount of two pebbles, \ie~$k=2$ for all edges. We will in some cases refer to this as \emph{weight-2} pebbling. To avoid trivialities, we assume that edges always have a weight of \emph{at least} two; $k\ge 2$.

\begin{definition}
Let $G$ be a graph.
A configuration $c$ on $G$ is a function $V(G)\to\Nat$.
\end{definition}

\begin{definition}
Let $G$ be a graph.
The pebbling step relation $c_1\to c_2$ between configurations $c_1$ and $c_2$ on $G$ holds if there is an edge $(u,v)$ such that $c_1(u)\ge k$ for $k=\weight_G(u,v)$, and $c_2$ results from $c_1$ after $k$ pebbles are removed from $u$ and one pebble is added to $v$:
$$
\forall w\in V(G):c_2(w)=\begin{cases}
c_1(u)-k&\If~w=u\\
c_1(v)+1&\If~w=v\\
c_1(w)&\Else
\end{cases}
$$
\end{definition}

We use an asterisk to denote the reflexive-transitive closure of this relation: $c_1\rtcto c_n$ means that there is a sequence of zero or more pebbling steps $c_1\to c_2\to\dots\to c_n$ from configuration $c_1$ to $c_n$. If, starting from a configuration $c$, it is possible to put $n$ pebbles on a target vertex~$t$ via pebbling steps, then we say that $c$ is $n$-fold $t$-solvable, or that $n$ pebbles can be \emph{moved}\footnote{Here `moving' $n$ pebbles also describes cases where the initial configuration already puts some number of pebbles on $t$, such that less or no additional pebbles need to be added to reach $n$.} to $t$. When $n=1$ we just write ``$t$-solvable''.

\begin{definition}
Let $G$ be a graph.
A configuration $c$ on $G$ is $n$-fold $t$-solvable if there is a configuration $c^\star$ such that $c\rtcto c^\star$ and $c^\star(t)\ge n$. As a first-order formula:
$$\solv{n}{t}{c}\defeq\Exists{c^\star}{c\rtcto c^\star\land c^\star(t)\ge n}$$
\end{definition}

\pagebreak
\begin{definition}
Let $G$ be a graph and let $t\in V(G)$.
The $n$-fold $t$-pebbling number $\npebb{n}{G,t}$ is the smallest number $p$ such that every configuration $c$ with $\size{c}\ge p$ is $n$-fold $t$-solvable. The pebbling number $\pebb{G}$ of $G$ is defined as the largest 1-fold pebbling number among its vertices.
\begin{align*}
\npebb{n}{G,t}&\defeq\min\Set{p\in\Nat}{\Forall{c}{\size{c}\ge p\implies\solv{n}{t}{c}}}\\
\pebb{G,t}&\defeq\npebb{1}{G,t}\\
\pebb{G}&\defeq\max\Set{\pebb{G,v}}{v\in V(G)}
\end{align*}
\end{definition}

\noindent
Let's use these definitions to prove some simple theorems.

\begin{theorem}
For a simple graph $G$ holds $\pebb{G^{(2)}}\ge\Count{V(G)}$.
\end{theorem}
\begin{proof}
Note that the weight of all edges in $G^{(2)}$ is 2 by definition.
Pick any vertex $t$ of $G$ and put a single pebble on each vertex except $t$. The resulting configuration contains $\Count{V(G)}-1$ pebbles and is not $t$-solvable since no pebbling steps can be applied.
\end{proof}

\begin{theorem}
\label{thm:pebb_K}
For $n\ge 1$ and $k\ge 1$ holds $\pebb{\completegraph{n}{k}}=(n-1)(k-1)+1$.
\end{theorem}
\begin{proof}
Note that all vertices of $\completegraph{n}{k}$ have the same pebbling number. Let $t$ be any target vertex. Placing $k-1$ pebbles at every vertex except $t$ results in a configuration that is not $t$-solvable, hence the pebbling number is larger than $(n-1)(k-1)$. If $(n-1)(k-1)+1$ pebbles are distributed over the vertices such that $t$ has 0 pebbles (otherwise the solution is trivial), then by the pigeonhole principle one vertex must have at least $k$ pebbles, such that one pebble can be moved to $t$.
\end{proof}

\subsection*{Trees}
Now we turn to something a bit more impressive: a formula for the $n$-fold pebbling number of a vertex in a tree. This formula was proven by Chung, and uses a so called \emph{maximum path-partition}.

\begin{definition}
Let~\,$T$ be a tree and let $r$ be a vertex of~~$T$. Let~~$T^\star$ be a directed graph obtained from $T$ by pointing all edges towards $r$. A path-partition rooted at $r$ is a sequence of disjoint vertex sets $V_1,V_2,\dots,V_n$ with a nonincreasing size, that each induce a directed path in~$T^\star$, such that ${V_1\cup V_2\cup\dots\cup V_n= V(T)\setminus\{r\}}$. A path-partition is maximum if the sequence of path sizes, $\Count{V_1},\Count{V_2},\dots,\Count{V_n}$, majorizes that of any other path-partition (rooted at the same vertex) lexicographically.
\end{definition}

\begin{theorem}
\cite{Chung1989}
Let $n\ge 1$ and $k\ge 2$. Let $T$ be a tree and let $r$ be a vertex of~~$T$.
If~~$V_1,V_2,\dots,V_m$ is a maximum path-partition of~~$T$ rooted at $r$, then:
$$\npebb{n}{T^{(k)},r}=n\cdot k^{\SmallCount{V_1}}+k^{\SmallCount{V_2}}+\dots+k^{\SmallCount{V_m}}-m+1$$
\end{theorem}
\pagebreak

\begin{figure}[h]
\includesvg{fig/example/tree_partition.svg}
\caption{Maximum path-partition}
\end{figure}

\begin{proof}
We prove this by induction on the number of vertices of~\,$T$. Suppose the theorem holds for all trees with fewer vertices than $T$. Remove $r$ from $T$ and obtain the subtrees $T_1,\dots,T_s$ that $T$ splits into, where $T_1$ contains the longest path directed towards $r$. For each $i\le s$,\linebreak let $v_i$ be the vertex of $T_i$ that is a neighbor of $r$. A maximum path-partition of~\,$T$ can be obtained from the union of, for every $i\le s$, a maximum path-partition of $T_i$ rooted at $v_i$ with $v_i$ added to its longest path.

Moving $n$ pebbles to $r$ can be achieved by any combination of moving $n_i$ pebbles to $v_i$ such that $\lfloor n_1/k\rfloor+\dots+\lfloor n_s/k\rfloor=n$. Note that in general, $\pi_{n+1}(G,v)-1$ is the largest number of pebbles that can be distributed over $G$ in such a way that it is impossible to move more than $n$ pebbles to $v$. Using this insight we can see that every sequence $n_1,\dots,n_s$ such that $\lfloor n_1/k\rfloor+\dots+\lfloor n_s/k\rfloor<n$ yields a lower bound for $\npebb{n}{T,r}$:

$$\npebb{n}{T,r}>\sum_{i\le s}\left(\npebb{n_i+1}{T_i,v_i}-1\right)$$

Let $c$ be a configuration on $T$ that is not $n$-fold $r$-solvable. Obtain a configuration $c^\ast$ by removing all pebbles on $r$ and adding them one of its neighbors. Note that $c^\ast$ is also not $n$-fold $r$-solvable. Split $c^\ast$ into $c_1,c_2,\dots,c_s$, where $c_i$ is a configuration on $T_i$, and for every $i\le s$, let $n_i$ be the maximum number of pebbles that can be moved to $v_i$ using the pebbles from $c_i$.
Note that $\size{c}=\size{c^\ast}\le \sum_{i\le s}\left(\npebb{n_i+1}{T_i,v_i}-1\right)$ and $\lfloor n_1/k\rfloor+\dots+\lfloor n_s/k\rfloor<n$.
We conclude from this that the maximum of the above lower bounds is strict:

$$\npebb{n}{T,r}-1=\max\Set{\sum_{i\le s}\left(\npebb{n_i+1}{T_i,v_i}-1\right)}{\sum_{i\le s}\lfloor n_i/k\rfloor<n}$$

Observe that $k^a+k^b\ge k^{a-1}+k^{b+1}$ if $a>b$. By the induction hypothesis, each\linebreak $\npebb{n_i+1}{T_i,v_i}$ can be replaced by an expression of the form $(n_i+1)k^x+c$, where $c$ is a constant and $x$ is the size of the longest path in $T_i$ (excluding $v_i$). Since $T_1$ contains the longest path, the maximum is achieved when $n_1=nk-1$ and $n_i=k-1$ for $1<i\le s$. Using some extra work we can now rewrite the above expression into the desired form.
\end{proof}

\pagebreak
\subsection*{Diameter-2 Graphs}
In a diameter-2 graph, any two vertices are either adjacent or have a common neighbor. The weight-2 pebbling number of such graphs was determined in \cite{Pachter1995} to be bounded by the number of vertices plus one. We present a simplified version of this proof. In the following, when it is written that \emph{we may assume $X$}, it is implied that the case \emph{not $X$} can easily be seen to be solvable. Note that $X$ may actually yield a contradiction, in which case \emph{not $X$} and the associated solution follow. In Figure \ref{fig:diameter2}, the gray area represents a part of the graph that is repeated, indexed by $i$.

\begin{theorem}
\cite{Pachter1995}
For a simple diameter-2 graph $G$:
$$\pebb{G^{(2)}}\le\Count{V(G)}+1$$
\end{theorem}
\begin{proof}
Let $c$ be a configuration on $G$ such that $\size{c}=\Count{V(G)}+1$, and let $t$ be a target vertex. We show that $c$ is $t$-solvable. Label the vertices with more than one pebble $x_1,x_2,\dots,x_m$. We may assume no $x_i$ is equal to $t$ or a neighbor of $t$. For each $x_i$, pick a neighbor $y_i$ that connects to $t$. We may assume that each $y_i$ has zero pebbles, and that all $y_i$ are distinct.

\begin{Figure}{fig/proof/diameter2_simplified.svg}
\caption{}
\label{fig:diameter2}
\end{Figure}

A simple calculation now implies that at least one vertex has 3 pebbles. We call this vertex $v$ and remove it from the $x_i$ vertices. We may assume that none of the $x_i$ vertices are adjacent to $v$. For each $x_i$, pick a neighbor $z_i$ that connects to $v$. We may assume that each $z_i$ has zero pebbles, that all $z_i$ are distinct, and that all $z_i$ are distinct from all $y_i$. A simple calculation now implies that at least one vertex has four pebbles, such that the configuration is $t$-solvable.
\end{proof}

Since the pebbling number is always at least the number of vertices, it follows from the previous result that, for every simple diameter-2 graph $G$, we have:
$$\pebb{G^{(2)}}\in\{\Count{V(G)},\Count{V(G)}+1\}$$
It is possible to characterize exactly which diameter-2 graphs have a pebbling number equal to their number of vertices. This was done in \cite{Clarke1997}, where the terminology \emph{Class-0} (for those simple graphs with a pebbling number equal to their number of vertices) and \emph{Class-1} (for the other diameter-2 graphs) originated.

\pagebreak
\section{Questions}
\label{sec:questions}
When Chung introduced pebbling numbers in her paper about hypercubes and the \Erdos-Lemke conjecture, she also asked a number of general questions about them. It is worth revisiting those questions, since they are the primer to many later publications. Chung's definition of graph pebbling is analogous to ours when ignoring edge direction, \ie making all edges bidirectional, and setting a standard edge weight of 2. Therefore, in the following questions, you should assume that all graphs are of this form.

\begin{question}
Is it true that $\pebb{G}=\max\{2^{\diam(G)},\Count{V(G)}\}$?
\end{question}

\emph{No.} The star graph $S_3$ with three leaves has 4 vertices and a diameter of 2, but $\pebb{S_3}=5$. The 7-cycle $C_7$ has 7 vertices and a diameter of 3, but $\pebb{C_7}=11$.

\begin{question}
Is it true that $\pebb{G\gprod H}=\pebb{G}\pebb{H}$?
\end{question}

\emph{No.}
Consider $C_3\gprod P_3$. We have $\pebb{C_3}=3$ and $\pebb{P_3}=4$, while $\pebb{C_3\gprod P_3}=9$. At the end of this chapter we will briefly explain a method based on CTL model checking that we implemented to automatically compute the pebbling number of $C_3\gprod P_3$.

\begin{question}
\emph{(Ronald Graham)}
Is is true that $\pebb{G\gprod H}\le\pebb{G}\pebb{H}$?
\end{question}

\emph{Unsolved.}
This elusive question is often referred to as \emph{Graham's conjecture}, after Ronald Graham who Chung attributes it to. A lot of graph pebbling research has been focused on this conjecture. The \emph{2-pebbling property}, which was suggested by Chung based on her work about the pebbling number of the hypercube, has yielded quite some results.

\begin{definition}
\label{def:2PP}
Let $G$ be a graph. $G$ has the 2-pebbling property (2PP) if for every configuration on $G$ with at least $2\pebb{G}-q+1$ pebbles, where $q$ equals the number of vertices that have at least one pebble, it is possible to move two pebbles to any target vertex.
\end{definition}

On many graphs, configurations that are more `spread out' need, in the worst case, fewer pebbles to be solvable. Several authors have used the 2-pebbling property to prove special cases of Graham's conjecture. An overview of some major results is given below.

\begin{table}[h]
\centering
\begin{tabular}{lll}
\toprule
The Cartesian product of \dots & Year & Reference\\
\midrule
\dots~$K_n$ with a 2PP graph & 1989 & \cite{Chung1989}\\
\dots~a sequence of trees & 1992 & \cite{Moews1992}\\
\dots~two 5-cycles & 1998 & \cite{Herscovici1998}\\
\dots~a tree with a 2PP graph & 2000 & \cite{Snevily2000}\\
\dots~a sequence of cycles (except $C_5$) and a 2PP graph & 2008 & \cite{Herscovici2008}\\
\dots~the Lemke graph\footnote{See Figure \ref{fig:L_labeled} on the next page.} with a tree or with $K_n$ & 2017 & \cite{Gao2017}\\
\bottomrule
\end{tabular}
\caption{Confirmed cases of Graham's conjecture.}
\label{tab:products}
\end{table}

\pagebreak
\begin{question}
Is it true that every graph has the 2-pebbling property?
\end{question}

\emph{No.}
Paul Lemke identified the smallest graph that does not have the 2-pebbling property, which is now called the \emph{Lemke graph} and is denoted as $L$. Two representations of the Lemke graph are shown below: the left depiction is standard in the literature, and the right depiction is the one we prefer to use. It can be determined that $\pebb{L}=8$. The configuration given in Figure \ref{fig:L_not_2PP} demonstrates that $L$ does not satisfy the 2-pebbling property.

\begin{Figure}{fig/graph/L_labeled.svg}
\caption{Depictions of the Lemke graph}
\label{fig:L_labeled}
\end{Figure}

\begin{Figure}{fig/example/L_not_2PP.svg}
\caption{$L$ does not satisfy the 2PP.}
\label{fig:L_not_2PP}
\end{Figure}

The Lemke graph is interesting because it challenges the methods used to find the pebbling number of graph products. The pebbling number of $L\gprod L$ is still unknown; if there exists a configuration of 64 pebbles on this graph that is unsolvable for some target vertex, then Graham's conjecture is false! Quite recently, with the help of linear optimization and so called \emph{weight functions}, it was determined that $\pebb{L\gprod L}\le 85$~\cite{Kenter2020}. In Chapter \ref{chapter:weight} we explain the basic ideas underlying these weight functions.


\begin{question}
Is $\pebb{C_5\gprod C_5\gprod\dots\gprod C_5}=5^n$,
where we take the product of $n$ 5-cycles?
\end{question}

\emph{Unsolved.}
So far this equality has only been proven for $C_5$ and $C_5\gprod C_5$ \cite{Herscovici1998}. Among cycles, the 5-cycle is a curious exception, because Graham's conjecture has been confirmed for products of cycles of \emph{any} other size (see Table \ref{tab:products}).

\pagebreak
\section{Computational Methods}
Is it possible to determine pebbling numbers of moderately sized graphs using a computer? A brute-force approach is not going to work; even on relatively small graphs the number of possible configurations is very high. Let's ignore automorphisms for a moment, and calculate the total number of ways in which $k$ pebbles can be distributed over $n$ distinct vertices:

\begin{theorem}
Let $G$ be a graph with $n$ vertices. For $k\in\Nat$:
$$\#{\Set{c:V(G)\to\Nat}{\size{c}=k}}={k+n-1\choose n-1}$$
\end{theorem}
\begin{proof}
Let $v_1,v_2,\dots,v_n$ be the vertices of $G$. We assign $k$ pebbles to these vertices. Pick numbers $p_1,p_2,\dots,p_{n-1}$ from $\{1,2,\dots,k+n-1\}$ such that $p_1<p_2<\dots<p_{n-1}$, this can be done in ${k+n-1\choose n-1}$ ways. Define a configuration $c$ as follows:
$$c(v_i)\defeq\begin{cases}
p_1-1&\If~i=1\\
p_i-(p_{i-1}+1)&\If~1<i<n\\
(k+n)-(p_{i-1}+1)&\If~i=n
\end{cases}$$
Now $\size{c}=k$. Note that every possible configuration of size $k$ can be picked in this way.
\end{proof}

Two configurations that only differ by an automorphism permuting the vertices admit the same solutions. We could divide the above formula by the number of automorphisms to get a rough estimate of the number of \emph{different} configurations. By this metric there are \emph{at least} 770 different configurations of 10 pebbles on the Petersen graph, and more than $3\times 10^{11}$ different configurations of 25 pebbles on $C_5\gprod C_5$\footnote{The graph $C_5\times C_5$ has 200 automorphisms.}. To date, no general algorithm for pebbling numbers has been found that can handle $C_5\gprod C_5$. It is worth noting that just solving individual configurations is already NP-complete~\cite{Milans2006}. At the end of Chapter~\ref{chapter:flow} we describe a general method to find pebbling solutions using integer programming.

\paragraph{Barely Sufficient Configurations}
In \cite{Sieben2010} a method is developed to find the set of \emph{barely sufficient} configurations using a recursive algorithm that applies inverse pebbling steps. A configuration is barely sufficient when it is solvable for some vertex, and it is no longer solvable when one pebble is removed from any vertex. The weight-2 pebbling number of all simple graphs with fewer than 10 vertices was computed using this method (combined with a graph simplification technique) in about a day. This publication is also the only one we found that describes graph pebbling using weighted graphs.

\paragraph{CTL model checking}
We implemented a method to compute pebbling numbers of small graphs using \emph{CTL model checking} and the NuSMV software. A CTL model describes, over a finite space of variables; a set of initial states, a state transition relation, and a specification about reachable states. The specification is verified by the model checker using a highly optimized algorithm based on binary decision diagrams. Both the set of initial states and the transition relation can be expressed as a boolean formula over the model variables.

\pagebreak
\begin{figure}[h]
\begin{tcolorbox}
\small\begin{verbatim}
MODULE main
DEFINE n := 4;
VAR c : array 1..3 of 0..n;
INIT c[1] + c[2] + c[3] = n

TRANS
 ( c[1]>1 & next(c[1])=c[1]-2 & next(c[2])=c[2]+1 & next(c[3])=c[3]   ) |
 ( c[2]>1 & next(c[1])=c[1]+1 & next(c[2])=c[2]-2 & next(c[3])=c[3]   ) |
 ( c[2]>1 & next(c[1])=c[1]   & next(c[2])=c[2]-2 & next(c[3])=c[3]+1 ) |
 ( c[3]>1 & next(c[1])=c[1]   & next(c[2])=c[2]+1 & next(c[3])=c[3]-2 ) |
          ( next(c[1])=c[1]   & next(c[2])=c[2]   & next(c[3])=c[3]   )

SPEC EF c[1] > 0
SPEC EF c[2] > 0
SPEC EF c[3] > 0
\end{verbatim}
\end{tcolorbox}
\caption{NuSMV model for $P_3$}
\end{figure}

\witregel

The above figure shows a NuSMV model for weight-2 pebbling on $P_3$. We want to check that starting at any configuration of four pebbles on $P_3$, it is possible to move one pebble to any target vertex. The variables \verb|c[1]|, \verb|c[2]| and \verb|c[3]| express the number of pebbles on each vertex. Since there can at most be four pebbles on any vertex, the value of these variables ranges from zero to four. The \verb|INIT| formula describing the initial set of states is true for any configuration of four pebbles. The \verb|TRANS| formula describing the transition relation between the current state, described by \verb|c[1]|, \verb|c[2]| and \verb|c[3]|, and the next state, described by \verb|next(c[1])|, \verb|next(c[2])| and \verb|next(c[3])|, is true for all valid pebbling steps and, for technical reasons, a transition where no pebbles are moved. For each vertex, a formula is defined using \verb|SPEC| to specify that one pebble can be moved to that vertex from any initial state. The specification formula \verb|EF|~$\varphi$ is true if, for all initial states, there is a sequence of transitions such that the formula $\varphi$ is true in the final state\footnote{See page 37 of \cite{NuSMVManual}.}.

\witregel

Figure \ref{fig:smv_L} on the next page shows a part of a NuSMV model for the Lemke graph. Here the initial states contain all configurations from which, to satisfy the 2-pebbling property, it should be possible to move two pebbles to any target vertex. When verifying this model, NuSMV finds a counterexample to $\verb|SPEC EF c[1] > 0|$, such as the one shown in Figure \ref{fig:L_not_2PP}. It takes NuSMV a few seconds to verify the weight-2 pebbling number of graphs like $P_3\gprod C_3$ and $L$, and doing so for the Petersen graph already takes around 10 seconds. Larger graphs quickly become intractable.

\begin{figure}[H]
\begin{tcolorbox}
\small\begin{verbatim}
MODULE main
DEFINE n := 8; p := 8;
VAR c : array 1..n of 0..2*p;

INIT
 c[1] + c[2] + c[3] + c[4] + c[5] + c[6] + c[7] + c[8] = 2*p + 1 -
 count(c[1]>0, c[2]>0, c[3]>0, c[4]>0, c[5]>0, c[6]>0, c[7]>0, c[8]>0)

TRANS
 ( c[1]>1 & next(c[1])=c[1]-2 & next(c[2])=c[2]+1 & ... ) |
 ( c[1]>1 & ... ) |
 ...

SPEC EF c[1] > 1
SPEC EF c[2] > 1
SPEC EF c[3] > 1
...
\end{verbatim}
\end{tcolorbox}
\caption{NuSMV model to check if $L$ satisfies the 2PP.}
\label{fig:smv_L}
\end{figure}

\begin{figure}[H]
\begin{tcolorbox}
\small\begin{verbatim}
-- specification EF c[1] > 1  is false
-- as demonstrated by the following execution sequence
Trace Description: CTL Counterexample
Trace Type: Counterexample
  -> State: 1.1 <-
    c[1] = 0
    c[2] = 0
    c[3] = 0
    c[4] = 1
    c[5] = 1
    c[6] = 1
    c[7] = 1
    c[8] = 8
    p = 8
    n = 8
-- specification EF c[2] > 1  is true
-- specification EF c[3] > 1  is true
...
\end{verbatim}
\end{tcolorbox}
\caption{NuSMV output}
\end{figure}

 \chapter{Weight Functions}
\label{chapter:weight}

In this chapter we use weight functions that assign weights to the \emph{vertices} of a graph. These weight functions should be distinguished from the edge-weight functions we introduced earlier. Weight functions enable us to recognize solvable pebbling configurations with linear inequalities, offering a powerful method for finding pebbling numbers. The use of weight functions and linear optimization was proposed by Hurlbert~\cite{Hurlbert2011}\cite{Hurlbert2017}. We only discuss weight functions applied to weight-2 pebbling, so for the duration of this chapter you should assume that the edge weight is always 2. We start by proving a key theorem, and showing how the pebbling number of various graphs can be determined using a set of weight functions.

\begin{definition}
Let $G$ be a graph and $t$ one of its vertices.
A function $w:V(G)\to\Frac$ is a weight function for $t$ if $w(t)=0$, and for every vertex $u\in\supp{w}$ such that $(u,t)\notin E(G)$, there is an edge $(u,v)$ such that $w(v)\ge2\cdot w(u)$.\footnote{In accordance with the literature, weight values are non-negative rationals. This can be useful since the neighbors of the target vertex could just as well be given a weight of 1, their neighbors a weight of \textonehalf, and so on.}
\end{definition}

Given a weight function, we can compute the weight of a configuration by multiplying, for each vertex, the number of pebbles with the weight assigned by the weight function, and adding the results. For a weight function $w$ and a configuration $c$ this is denoted as $\size{w\fprod c}$. The defining property of a weight function makes sure that pebbles on a positively weighted vertex not adjacent to the target vertex can be moved to a neighboring vertex (via pebbling steps) without the configuration as a whole losing weight. This way we can prove the following result, which Hurlbert calls the \emph{Weight Function Lemma}.

\begin{lemma}
\label{lem:weight_functions}
Let $G$ be a graph and $t$ one of its vertices.
If $w$ is a weight function for $t$, and $c$ is a configuration on $G$ such that $\size{w\fprod c}>\size{w}$, then $c$ is $t$-solvable.
\end{lemma}
\begin{proof}
Let $c$ be a configuration on $G$ such that $\size{w\fprod c}>\size{w}$. Determine a vertex $u$ with $c(u)\ge 2$. If $(u,t)\in E(G)$, then we can move one pebble to $t$, and so $c$ is $t$-solvable. Otherwise, determine an edge $(u,v)$ such that $w(v)\ge2\cdot w(u)$, and apply one pebbling step to this edge, removing two pebbles from $u$ and adding one to $v$. For the updated configuration $c'$ we have $\size{w\fprod c'}\ge\size{w\fprod c}>\size{w}$. Repeat this procedure until a pebble has been moved to $t$. This should only take a finite number of steps, since there is only a finite number of pebbles. We conclude that $c$ is $t$-solvable.
\end{proof}

We are going to use sums of weight functions to determine the pebbling number of several graphs, including the Petersen graph and cycle graphs. Suppose a sequence $w_1,w_2,\dots,w_n$ of weight functions for a vertex $t$ in a graph $G$ is given. If the sum $w=w_1+w_2+\dots+w_n$ is positive for all vertices except $t$, then we say that $w$ \emph{covers} $G$, and if all positive values of $w$ are equal, then we say that $w$ covers~$G$ \emph{uniformly}. The smallest positive value of $w$ is called its \emph{minimum weight}.

\pagebreak
\begin{lemma}
\label{lem:weight_fn_sum}
Let $G$ be a graph and $t$ one of its vertices.
Let $w_1,w_2,\dots,w_n$ be weight functions for $t$, and let $w=w_1+w_2+\dots+w_n$. If $w$ covers $G$ and $m$ is its minimum weight, then:
$$\pebb{G,t}\le\lfloor\size{w}/m\rfloor + 1$$
\end{lemma}
\begin{proof}
Let $c$ be a configuration on $G$ with $\size{c}\ge\lfloor\size{w}/m\rfloor + 1$. We prove that $c$ is $t$-solvable. Suppose $t$ has no pebbles, \ie~$c(t)=0$. Now $\size{w\fprod c}\ge\size{c}\cdot m>\size{w}$, which we can expand into $\size{w_1\fprod c}+\dots+\size{w_n\fprod c}>\size{w_1}+\dots+\size{w_n}$.
Determine $i\le n$ such that $\size{w_i\fprod c}>\size{w_i}$ and apply Lemma \ref{lem:weight_functions} to conclude that $c$ is indeed $t$-solvable.
\end{proof}

\begin{theorem}
For the Petersen graph $P$ holds:
$$\pebb{P}=10$$
\end{theorem}
\begin{proof}
Since $P$ has 10 vertices, $\pebb{P}\ge 10$. Note that all vertices of $P$ have the same pebbling number. Pick some target vertex $t$. For each neighbor of $t$ there is a weight function like the one shown below. Adding these three weight functions together covers $P$ uniformly with fours. Using Lemma \ref{lem:weight_fn_sum} we find $\pebb{P,t}\le\lfloor 36/4\rfloor+1=10$.
\begin{center}
\vspace{2mm}
\includesvg{fig/proof/weight_fn_Petersen.svg}
\end{center}
\end{proof}

\begin{theorem}
For the complete bipartite graph $K_{m,n}$ with $m,n\ge 2$ holds:
$$\pebb{K_{m,n}}=m+n$$
\end{theorem}
\begin{proof}
The sum of the two weight functions shown below covers the graph uniformly with twos. Using Lemma \ref{lem:weight_fn_sum} we find $\pebb{K_{m,n}}\le m+n$.
\begin{center}
\vspace{2mm}
\includesvg{fig/proof/weight_fn_bipartite.svg}
\end{center}
\end{proof}

\pagebreak
\begin{theorem}
For the Lemke graph $L$ with $v_3$ as in Figure \ref{fig:L_labeled} holds:
$$\pebb{L,v_3}=8$$
\end{theorem}
\begin{proof}
Use the three weight functions shown below.
\begin{center}
\vspace{2mm}
\includesvg{fig/proof/weight_fn_Lemke.svg}
\end{center}
\end{proof}

\begin{theorem}
For the even cycle $C_{2n}$ with $n\ge 1$ holds:
$$\pebb{C_{2n}}=2^n$$
\end{theorem}
\begin{proof}
Make two mirror copies of the weight function shown below, and call their sum $w$. The minimal weight of $w$ is 2 and $\size{w}=2(2^n-1)$. Using Lemma \ref{lem:weight_fn_sum} we find $\pebb{C_{2n}}\le 2^n$. This is strict, since the configuration that puts $2^n-1$ pebbles on the vertex furthest from the target is unsolvable.
\begin{center}
\vspace{2mm}
\includesvg{fig/proof/weight_fn_even_cycle.svg}
\end{center}
\end{proof}

\begin{theorem}
For the odd cycle $C_{2n+1}$ with $n\ge 1$ holds:
$$\pebb{C_{2n+1}}=2\lfloor{2^{n+1}/3}\rfloor+1$$
\end{theorem}
\begin{proof}
Make two mirror copies of the weight function shown below, and call their sum~$w$. The minimal weight of $w$ is 3 and $\size{w}=2(2^{n+1}-1)$. It follows from Lemma \ref{lem:weight_fn_sum} that $\pebb{C_{2n+1}}\le\lfloor2(2^{n+1}/3-1/3)\rfloor+1=2\lfloor{2^{n+1}/3}\rfloor+1$, where the last equality can be determined by analyzing the possible remainders of $2^{n+1}/3$. To see that this bound is strict, put $\lfloor{2^{n+1}/3}\rfloor$ pebbles on each of the two vertices furthest from the target, and note that no solution is possible via either side of the cycle since $\lfloor{2^{n+1}/3}\rfloor+\lfloor{2^n/3}\rfloor<2^n$.
\begin{center}
\vspace{2mm}
\includesvg{fig/proof/weight_fn_odd_cycle.svg}
\end{center}
\end{proof}

\section{Limitations}
So far weight functions have been quite instrumental, but there are some serious limitations. For some configurations there is no weight function to solve it using the procedure in the proof of Lemma \ref{lem:weight_functions}. If for a configuration $c$ there is a weight function $w$ with $\size{w\fprod c}>\size{w}$, then we say that $c$ can be solved using a weight function, and if every configuration $c$ on a graph $G$ with $\size{c}=\pebb{G}$ can be solved using a weight function, then we say that $G$ can be solved using weight functions. In this section we look at several graphs that are \emph{not} solvable using weight functions.

\subsection*{Splitting Structures}
Consider the configurations on the cube and the Lemke graph in the figure below. Both are solvable, but in order to solve them the heap of 5 pebbles must be divided over two neighbors to take advantage of the extra pebbles which are placed there. There is no other solution. Hurlbert calls these situations \emph{splitting structures}.

\begin{Figure}{fig/example/splitting_structures.svg}
\caption{}
\label{fig:splitting}
\end{Figure}

In the procedure that we used to prove Lemma \ref{lem:weight_functions} it does not matter how pebbles are distributed over neighbors, as long as the configuration weight is not decreased. In both of the above configurations, it is impossible to come up with a weight function to solve it. Both configurations also have the same number of pebbles as the (weight-2) pebbling number of their respective graphs. If the pebbling number of either of these graphs could be determined using just weight functions, then there should also be a weight function to solve these specific configurations. Hence we are forced conclude that the pebbling number of neither the cube nor the Lemke graph can be found using weight functions alone.

\witregel

In \cite{Cranston2017} this issue is bypassed by introducing a special weight function called a `lollipop' that contains a path attached to an even cycle. On a lollipop, the vertex furthest from the target gets slightly more weight than it can have in an ordinary weight function. This breaks the simple proof of Lemma \ref{lem:weight_functions}, since the configuration weight \emph{decreases} when a pebbling step removes pebbles from this vertex.

\pagebreak
But in Lemma 3 of \cite{Cranston2017} the authors show that their lollipop weight function $w_l$ has the same property as ordinary weight functions: If $\size{w_l\fprod c}>\size{w_l}$, then $c$ is solvable. An example of a lollipop weight function on the cube, solving the problematic configuration from Figure~\ref{fig:splitting}, is shown below.

\begin{Figure}{fig/example/weight_fn_lollipop.svg}
\caption{A lollipop weight function}
\label{fig:lollipop}
\end{Figure}

\subsection*{Balanced Configurations}
Consider the two configurations on $S_3$ below; the left one is solvable and the right one is not. Since the left configuration is the average of two unsolvable configurations like the right one, it cannot be solved using a weight function.

\SimpleFigure{fig/example/balanced_S3.svg}{}

\noindent This example has fewer than $\pebb{S_3}$ pebbles, and it is still possible to show that $\pebb{S_3}=5$ using weight functions. This led us to ask the following question: \emph{Given a graph $G$, can we find $\pebb{G}$ using just weight functions if all configurations $c$ with $\size{c}=\pebb{G}$ can be solved without splitting?}

\witregel

It appears the answer is \emph{no}. It has been shown that $\pebb{C_5\gprod C_5}=25$, and as far as we know all configurations of 25 pebbles on $C_5\gprod C_5$ can be solved without splitting\footnote{A careful analysis of \cite{Herscovici1998} might be able to confirm this.}. But we did find a configuration of 25 pebbles that is not solvable using a weight function. The offending configuration, and a weight function that comes closest to solving it, are shown in Figure~\ref{fig:balanced}.

\begin{Figure}{fig/example/balanced_C5xC5.svg}
\vspace{1mm}
\caption{A configuration (left) and a weight function (right) on $C_5\gprod C_5$.}
\label{fig:balanced}
\end{Figure}

\section{Linear Programming}
Every weight function $w$ gives rise to a linear inequality of the form $\size{w\fprod c}\le\size{w}$ that holds if a configuration $c$ \emph{cannot} be solved with it. Linear programming software can efficiently determine the optimal solution of a linear objective function subject to linear constraints. When we use the configuration size as the objective function, and the linear inequalities belonging to a set of weight functions as constraints, then the maximal solution gives the size of the largest possible configuration for which none of the weight functions work; all larger configurations will violate one of the weight function constraints, and are therefore solvable.

\paragraph{Fractional Relaxation}
A pebbling configuration consists of integer amounts; we do not split up pebbles between vertices, and thus we need an LP solver that finds integer solutions. In linear programming, finding integer solutions is much more computationally intensive than finding fractional solutions, and therefore our problem becomes much easier to compute if we `relax' it to configurations with fractional pebble amounts. Since integer configurations are a subset of fractional configurations, this will still yield a useful bound.

\newcommand{\maxconf}{c_{\mathrm{max}}}
\begin{lemma}
Let $G$ be a graph and $t$ one of its vertices.
If $w_1,w_2,\dots,w_n$ are weight functions for $t$, and a fractional configuration $\maxconf:V(G)\to\Frac$ maximizes $\size{\maxconf}$ while satisfying $\size{w_i\fprod\maxconf}\le\size{w_i}$ for all $i\le n$, then:
$$\pebb{G,t}\le\lfloor\size{\maxconf}\rfloor+1$$
\end{lemma}
\begin{proof}
Let $c$ be a configuration on $G$ with $\size{c}\ge\lfloor\size{\maxconf}\rfloor+1$. Then $\size{c}>\size{\maxconf}$. Because $\maxconf$ is maximal, there must exist an $i\le n$ such that $\size{w_i\fprod c}>\size{w_i}$. It now follows from Lemma~\ref{lem:weight_functions} that $c$ is $t$-solvable.
\end{proof}

\pagebreak
\paragraph{Dual Solutions}
The \emph{dual solution} of an LP solver reveals the constraints by which the optimal corner of the feasibility polytope\footnote{The feasibility polytope is the region containing all points that satisfy the given linear constraints.} is bounded. It is possible to enter many weight functions into an LP solver, and then use the dual solution to find out which weight functions are sufficient to obtain the optimal bound. Hurlbert did this in \cite{Hurlbert2011}; he first picked a random graph on 15 vertices, then a program generated over 20 thousands weight functions for the optimization problem, and finally the CPLEX solver found an optimum based on just 11 of these weight functions. In this method, the feasibility polytope for a graph with $n$ vertices is $n$-dimensional, where each dimension represents the number of pebbles on one of its vertices. The corners of a 3-dimensional polytope are almost always at the intersection of 3 of its faces, and in general a linear optimization problem with $n$ variables typically requires at most $n$ constraints in its certificate.

\vspace{1cm}

\begin{figure}[H]
\includesvg[width=.6\textwidth]{fig/misc/simplex.svg}
\vspace{3mm}
\caption{Integer solutions}
\end{figure}

\chapter{Pebble Flows}
\label{chapter:flow}

What can we say about the direction in which pebbles are moved? Hurlbert calls solutions \emph{greedy} when pebbles are moved strictly closer to the target vertex. For some graphs, for example for cubes, there is a greedy solution for all configurations with a size equal to the graph's pebbling number. But sometimes pebbles actually have to be moved \emph{away} from the target vertex. This is the case in the weight-2 example shown below. It is possible to determine that the weight-2 pebbling number of this graph is 10 using weight functions.

\SimpleFigure{fig/example/detour_3_bridges.svg}{}

\emph{Do we ever need to move pebbles back and forth between two vertices?} This does not seem very useful; if pebbles that are moved to some vertex can be returned later, then how are they essential to a solution? In this chapter we show that, indeed, to solve any particular configuration it is never needed to move pebbles back-and-forth between any two vertices. To prove this, we regard solutions as \emph{flow networks}. We believe that this approach might also be useful for other graph pebbling problems.

\begin{definition}
Let $G$ be a graph.
A pebble flow $\fF$ on $G$ is a pair $(c,f)$ of a configuration $c$ on $G$, and a flow function $f:V(G)\times V(G)\to\Nat$. The value $f(u,v)$ represents the number of pebbles that are added to $v$ using pebbling steps via the edge $(u,v)$. If $(u,v)\notin E(G)$ then $f(u,v)=0$. The configuration of $\fF$ is also written as $c_\fF$ and the flow from $u$ to $v$ as $\fF(u,v)$.
\end{definition}

A pebble flow combines a configuration with a solution expressed in flow values. We can compute several useful values from this, including the total inflow and outflow of a vertex, the \emph{weighted} outflow of a vertex which is the number of pebbles that actually need to be removed to realize its outflow, and the \emph{excess} of a vertex which is the number of pebbles that are left after the inflow and outflow are satisfied. We will define these terms more precisely on the next page.

\pagebreak
\begin{definition}
Let $G$ be a graph.
For a pebble flow $\fF$ on $G$ and a vertex $v\in V(G)$, the inflow $\inflow{\fF}{v}$, outflow $\outflow{\fF}{v}$, weighted outflow $\outfloww{\fF}{v}$, and excess $\excess{\fF}{v}$ of $v$ are defined as:
\begin{align*}
\inflow{\fF}{v}&\defeq\sum_{u\in V(G)}\fF(u,v)\\
\outflow{\fF}{v}&\defeq\sum_{u\in V(G)}\fF(v,u)\\
\outfloww{\fF}{v}&\defeq\sum_{u\in V(G)}\weight_G(v,u)\cdot\fF(v,u)\\
\excess{\fF}{v}&\defeq c_\fF(v)+\inflow{\fF}{v}-\outfloww{\fF}{v}
\end{align*}
\end{definition}

\noindent
We extend the step relation for pebbling configurations to pebble flows as follows:

\begin{definition}
Let $G$ be a graph.
The pebble flow step relation $\fF_1\to \fF_2$ between two pebble flows $\fF_1$ and $\fF_2$ on $G$ holds if there is an edge $(u,v)$ such that $\fF_1(u,v)\ge 1$ and $c_{\fF_1}(u)\ge k$ for $k=\weight_G(u,v)$, and:
\begin{align*}
\forall w:c_{\fF_2}(w)&=\begin{cases}
c_{\fF_1}(u)-k&\If~w=u\\
c_{\fF_1}(v)+1&\If~w=v\\
c_{\fF_1}(w)&\Else
\end{cases}\\[2mm]
\forall x,y:\fF_2(x,y)&=\begin{cases}
\fF_1(u,v)-1&\If~(x,y)=(u,v)\\
\fF_1(x,y)&\Else
\end{cases}
\end{align*}
\end{definition}

\section{Properties}
To translate a step-by-step solution into a pebble flow, you have to count the number of pebbling steps along every edge. For every vertex, there will never flow more pebbles out (the weighted outflow) than the number of pebbles put there by the initial configuration plus the number of pebbles flowing in. This means that the excess of every vertex is non-negative. We call this property \emph{feasible}. Clearly, a solvable configuration corresponds to a feasible pebble flow.

\begin{definition}
A pebble flow $\fF$ is feasible if~~$\Forall{v}{\excess{\fF}{v}\ge 0}$.
\end{definition}

\begin{lemma}
Let $G$ be a graph and $t$ one of its vertices.
If $c$ is an $n$-fold $t$-solvable configuration, then there exists a feasible pebble flow~$\fF$ with $c_\fF=c$ and $\excess{\fF}{t}\ge n$.
\end{lemma}
\begin{proof}
Left to the reader.
\end{proof}

The above property is rather trivial, but can we turn it around? Does every feasible pebble flow correspond to a step-by-step solution? If we can prove that this is true, then we are very close to answering the question about back-and-forth movement between vertices.

\pagebreak
\subsection*{Realization}
We quickly run into the following problem: it is possible to create a kind of deadlock using a cycle. This is illustrated below, where each edge has a weight of 2 and a flow value of 1. While this pebble flow is feasible, it is \emph{not} possible to apply further pebbling steps and reduce all flow values to zero.

\SimpleFigure{fig/example/flow_cycle_triangle.svg}{Cyclic flow}

However, in a situation like this there isn't really anything to gain. The excess of all vertices is zero, and each vertex already has more than zero pebbles, so why not stop here? That is exactly what we are going to do: we will show that for any pebble flow it is possible to take steps forward until every vertex has at least as many pebbles as its excess, at which point we say that the pebble flow is \emph{realized}.

\begin{definition}
A pebble flow $\fF$ is realized if~~$\Forall{v}{c_\fF(v)\ge\excess{\fF}{v}}$.
\end{definition}

\begin{definition}
A pebble flow $\fF$ is realizable if $\fF\rtcto\fF^\star$ for some $\fF^\star$ that is realized.
\end{definition}

We found two methods to prove that all feasible pebble flows are realizable. The first method is to remove all flow cycles, leaving an acyclic pebble flow that can be realized starting at its leaves. Note that such a pebble flow will always have at least one leaf; a vertex with no inflow, and that, since it is feasible, this vertex has enough pebbles to fully realize its outflow.
The second method we found doesn't need to know anything about flow cycles; it just shows that, as long as the pebble flow is not fully realized, there must be a vertex somewhere with enough pebbles for a next pebbling step. Since there can only be a finite number of steps until all pebbles run out, a realized configuration should eventually be reached. We will work out the second method in more detail.

\begin{lemma}
\label{lem:flow_balance}
For a pebble flow $\fF$ on $G$, the following holds:
$$
\sum_{v\in V(G)}\inflow{\fF}{v}=
\sum_{v\in V(G)}\outflow{\fF}{v}
$$
\end{lemma}
\begin{proof}
Every flow value counts as inflow for one vertex and outflow for another vertex. Thus the sum of the inflow of all vertices is equal to the sum of the outflow of all vertices.
\end{proof}

\pagebreak
\begin{lemma}
\label{lem:flow_realizeable}
Every feasible pebble flow is realizable.
\end{lemma}
\begin{proof}
Suppose $\fF$ is a pebble flow on a graph $G$, and $\fF$ is not realized. Determine a vertex $v$ such that $c_\fF(v)<\excess{\fF}{v}$. Since the excess has not been realized, $v$ must still be waiting for more pebbles than it needs to send away: $\inflow{\fF}{v}>\outfloww{\fF}{v}\ge\outflow{\fF}{v}$. It follows from Lemma \ref{lem:flow_balance} that there also exists a vertex in $G$ where this inequality is reversed. Determine $w$ such that $\inflow{\fF}{w}<\outflow{\fF}{w}$. Let $k$ be the minimum weight of an outflow of $w$, such that: $$\outfloww{\fF}{w}\ge k\cdot\outflow{\fF}{w}\ge k\cdot(\inflow{\fF}{w}+1)\ge\inflow{\fF}{w}+k$$
It follows from the feasibility of $\fF$ that $c(w)\ge k$, and thus that $w$ has enough pebbles for one pebbling step. Repeat this process until the pebble flow is realized. Note that, as long as it is not realized, there must be enough pebbles somewhere to apply a pebbling step, and that we can only repeat this a finite number of times because every step reduces the total number of pebbles.
\end{proof}

\subsection*{Unidirectional}
We have shown that pebble flows correspond to step-by-step solutions, and step-by-step solutions correspond to pebble flows. To answer our question about back-and-forth movement it only remains to see that any solution gives rise to a solution without back-and-forth movement. We call such a solution \emph{unidirectional}.

\begin{theorem}
Let $G$ be a graph. If a configuration $c$ is $n$-fold $t$-solvable, then there is a unidirectional solution.
\end{theorem}
\begin{proof}
Suppose $c$ is $n$-fold $t$-solvable. Determine a feasible pebble flow $\fF$ corresponding to a solution of $c$, such that $\excess{\fF}{t}\ge n$. We define a new pebble flow $\fF^\star$ with the same configuration that only describes the net flow between any two vertices:
$$\fF^\star(u,v)\defeq\fF(u,v)-\min\{\fF(u,v),\fF(v,u)\}$$
This pebble flow is feasible \emph{and} unidirectional. The excess does not decrease at any vertex, and in particular we still have $\excess{\fF^\star}{t}\ge n$. Using Lemma \ref{lem:flow_realizeable} we can find a step-by-step solution corresponding to $\fF^\star$ that puts at least $n$ pebbles on $t$. In this solution there is no back-and-forth movement.
\end{proof}

\pagebreak
\section{Other Applications}
Transforming a pebbling solution into a flow network removes the temporal aspect of steps, and allows the solution to be presented \emph{at once}. The feasibility property forms a conjunction of linear inequalities; an integer programming problem! This provides a straightforward method to solve pebbling configurations.

\paragraph{Optimal Pebbling}
The \emph{optimal} pebbling number of a graph is the minimum size of a configuration that is solvable for all of its vertices, \ie such that one pebble can be moved to any desired vertex. There are several unresolved questions about the optimal pebbling number. Curiously, the weight-2 optimal pebbling number of square grids is not yet known. A lower bound was established quite recently in \cite{Győri2020}:
$$\optpebb{\pathgraph{m}{2}\gprod\pathgraph{n}{2}}\ge {2\over 13}mn$$
An asymptotic upper bound is obtained in~\cite{Győri2017} using a configuration on an infinite grid that, on average, spreads every pebble over $3.5$ vertices (the \emph{covering ratio}). This configuration is shown in the figure below.

\vspace{5mm}
\SimpleFigure{fig/example/optimal_grid.svg}{Configuration with covering ratio 3.5.}

\pagebreak
\paragraph{Optimal configurations via integer programming}
By expressing solutions as pebble flows, finding an optimal configuration can be achieved with integer programming. We tried using the Z3 solver\footnote{\url{https://www.microsoft.com/en-us/research/project/z3-3/}} to find optimal configurations for small grids, but this already appeared to be intractable\footnote{Computing a solution for a $3\times 4$ grid using the Z3 solver took about 10 seconds on an AMD~Ryzen~7~5700G~CPU. We ran a similar calculation for a $4\times 4$ grid, but we terminated it after about ten minutes.} for grids larger than $3\times 4$. Perhaps it is possible to use smarter solving heuristics to find optimal configurations for larger grids using this method. For example, the distance over which pebbles are moved across the grid could be restricted. This problem might also be an interesting benchmark for LP solvers.

\chapter{Hypercubes}
\label{chapter:hypercubes}

We now turn to the problem that initiated the study of graph pebbling: the pebbling number of the hypercube. Chung's method to prove the \Erdos-Lemke conjecture, which we stated in the introduction, uses the pebbling number of a hypercube with arbitrary edge weights under the restriction that parallel edges have the same weight. Given the weights $k_1,k_2,\dots,k_n$, this corresponds to the following $n$-dimensional hypercube:
$$\pathgraph{2}{k_1}\gprod\pathgraph{2}{k_2}\gprod\dots\gprod\pathgraph{2}{k_n}$$
Chung determined that the pebbling number of this graph is the product of the weights; $k_1\cdot k_2\cdots k_n$, but the proof she presents leaves some details to the reader, and we actually got stuck trying to work it out. We filled this gap by developing a slightly different proof, which is presented in this chapter.

\section{Analysis}
We start by analyzing Chung's proof in \cite{Chung1989} that an $n$-dimensional hypercube with a standard edge weight of 2 has pebbling number $2^n$. We have included a copy of this proof below in its original form, with some minor adjustments.

\begin{theorem}
\label{thm:hypercube_2}
\cite{Chung1989} In an $n$-cube with a specified vertex $v$, the following are true:
\begin{enumerate}[(i)]
\item If~~$2^n$ pebbles are assigned to vertices of the $n$-cube, one pebble can be moved to $v$.
\item Let $q$ be the number of vertices that are assigned an odd number of pebbles. If there are in total more than $2^{n+1}-q$ pebbles, then two pebbles can be moved to $v$.
\end{enumerate}
\end{theorem}
\begin{proof}
The proof is by induction on $n$. It is trivially true for $n=0$. Suppose it is true for~$n$, we will prove it for $n+1$. Partition the $(n+1)$-cube into two $n$-cubes called $M_1$ and $M_2$, where $v$ is in $M_1$. Let $v'$ denote the vertex in $M_2$ adjacent to $v$. The edges between $M_1$ and $M_2$ form a perfect matching. Suppose $M_i$ contains $p_i$ pebbles with $q_i$ vertices having an odd number of pebbles, for $i=1,2$.

\item\textbf{(i)}~
Suppose there are $p\ge 2^{n+1}$ pebbles assigned to vertices of the $(n+1)$-cube. If $p_1\ge 2^n$, then by induction one pebble can be moved to $v$ in $M_1$. We may assume $p_1<2^n$ and we consider the following two cases:

\item\textbf{Case (a1)}~~$q_2>p_1$.~
Since $p_2=p-p_1>2^{n+1}-q_2$, by induction from (ii) in $M_2$ two pebbles can be moved to $v'$. Therefore one pebble can be moved to $v$.

\item\textbf{Case (a2)}~~$q_2\le p_1$.~
We can move at least $(p_2-q_2)/2$ pebbles from $M_2$ to $M_1$. This results in a total of $p_1+(p_2-q_2)/2\ge p_1+(p_2-p_1)/2=(p_1+p_2)/2\ge 2^n$ pebbles in $M_1$. By induction, we can then move one pebble to $v$.

\pagebreak
\item\textbf{(ii)}~
Suppose there are $p=p_1+p_2>2^{n+2}-q_1-q_2$ pebbles assigned to vertices of the $(n+1)$-cube. We want to show that two pebbles can be moved to $v$. We consider three cases:

\item\textbf{Case (b1)}~~$p_1>2^{n+1}-q_1$.~
By induction from (ii) in $M_1$, two pebbles can be moved to $v$.

\item\textbf{Case (b2)}~~$2^n\le p_1\le 2^{n+1}-q_1$.~
Since $p_1\ge 2^n$, by induction from (i) in $M_1$ one pebble can be moved to $v$. Note that $p_2=p-p_1>2^{n+2}-q_1-q_2-p_1\ge 2^{n+1}-q_2$. Therefore, by induction from (ii) in $M_2$, two pebbles can be moved to $v'$, such that a second pebble can be moved to~$v$.

\item\textbf{Case (b3)}~~$p_1< 2^n$.~
For any integer $t$ satisfying $p_2\ge q_2+2t$, it is possible to move $t$ pebbles to $M_1$ while $p_2-2t$ pebbles remain in $M_2$. Note:
$$p_2>2^{n+2}-q_1-q_2-p_1={(2^{n+1}-q_2)}+(2^{n+1}-q_1-p_1)\ge q_2+(2^{n+1}-q_1-p_1)$$
The last inequality follows since $q_2$ is at most the number of vertices in $M_2$. Now taking $t$ to be $2^n-(p_1+q_1)/2$, we move $t$ pebbles to $M_1$, leaving more than $2^{n+1}-q_2$ pebbles in~$M_2$. In $M_1$ there are now $p_1+2^n-(p_1+q_1)/2=2^n+(p_1-q_1)/2\ge 2^n$ pebbles. We can then move one pebble to $v$ in $M_1$ and move two pebbles to $v'$ in $M_2$, enabling us to move a second pebble to $v$.
\end{proof}

\subsection*{Odd 2-Pebbling Property}
\label{sec:2-pebbling_property}
Where does property (ii) in Theorem \ref{thm:hypercube_2} come from? One way to look at it is like this; if in $M_1$ we want to move one pebble to $v$, then the first thing we could try is moving as many pebbles from $M_2$ to $M_1$ as we can, and seeing if this results in enough pebbles on $M_1$ to apply the induction hypothesis. The number of pebbles remaining in $M_2$ will be equal to $q_2$; the number of vertices that have an odd number of pebbles.
This is exactly what is happening in case (a2) of the proof. The above strategy works as long as $q_2$ does not exceed the total number of pebbles on $M_1$. If it does, then $M_2$ has enough pebbles to apply (ii), as is done in case (a1). Property (ii) is a kind of `residual case' resulting from the above strategy, and it is sometimes called the \emph{odd 2-pebbling property}.

\witregel

When generalizing the weight of the edges between $M_2$ and $M_1$, the number of pebbles that are left behind on $M_2$ in case (a2) is no longer equal to~$q_2$. If the weight between $M_2$ and $M_1$ is $k$, then in the worst case $k-1$ pebbles are left behind on \emph{every} vertex. We are going to `upgrade' the proof of Theorem~\ref{thm:hypercube_2} to hypercubes with varying weights using an alternative for property~(ii).

\subsection*{Pebble Extraction}
Could we prove Theorem~\ref{thm:hypercube_2} when (ii) is replaced with the normal 2-pebbling property? Going forward we will use $\suppsize{c}$ to denote $\Count{\supp{c}}$, \ie~the number of vertices on which the configuration~$c$ puts at least one pebble. Let $c$ be a configuration on the $n$-cube, can we prove the following, using the same case distinctions and strategies as the proof of Theorem~\ref{thm:hypercube_2}?

$$\text{If $\size{c}\ge 2^{n+1}-\suppsize{c}+1$, then $c$ is 2-fold $v$-solvable.}$$

\pagebreak
\noindent The answer is \emph{no}, as is demonstrated by the following example on the cube denoted by $Q_3$. Let $Q_3$ be composed of the squares $M_1$ and $M_2$ with corresponding configurations $c_1$ and $c_2$ as shown below. Note that $\pebb{C_4}=4$ and $\pebb{Q_3}=8$. The following holds:
$$|c_1|+|c_2|=10\ge2\pebb{Q_3}-\suppsize{c_1}-\suppsize{c_2}+1$$

\begin{figure}[H]
\includesvg{fig/example/extraction_counterexample_1.svg}
\caption{}
\end{figure}

According to the 2-pebbling property we should now be able to move two pebbles to any vertex in $M_1$. We are trying to follow the same steps as Chung's proof, so we apply case (b3) since $|c_1|=3<\pebb{C_4}$. According to (b3) we should first move one pebble to $M_1$ using two pebbles from $M_2$, and then apply the 2-pebbling property inductively in $M_2$. But removing two pebbles from a vertex of $M_2$ \emph{decreases} the support of $c_2$, and the resulting configuration no longer satisfies the condition needed to apply the 2-pebbling property. We have to do more to generalize the proof.

\witregel

Chung adds \emph{pebble extraction} to the 2-pebbling property to solve this: If a configuration $c$ on the $n$-cube contains at least $2^{n+1}-\suppsize{c}+1+m$ pebbles, then it is possible to move two pebbles to a target vertex \emph{and also} move $\lfloor m/k\rfloor$ pebbles\footnote{Theorem 3 in \cite{Chung1989} actually allows a combination of pebbling steps with different weights.} to an adjacent $n$-cube using pebbling steps of weight $k$. Note that the latter is \emph{not} possible with every configuration of just $m$ pebbles, so there must be some interaction between these two parts.

We encountered a restriction on the weight of these extra pebbling steps that Chung does not mention. She proves that it works for $k=2$ if all edges have a weight of 2, but when generalizing $k$ we run into a problem. Let $c_3$ be the configuration below. For $m=4$ we have:
$$\size{c_3}\ge 2\pebb{C_4}-\suppsize{c_3}+1+m$$

\SimpleFigure{fig/example/extraction_counterexample_2.svg}{}

If we pick $k=4$ then $\lfloor m/k\rfloor=1$, and according to the above extraction property we should be able to find a vertex from which we can remove 4 pebbles. But $c_3$ does not have such a vertex, so it is not possible to move one pebble to an adjacent $n$-cube connected by edges of weight 4. To solve this problem we will apply induction from the dimension with the largest to the dimension with the smallest weight.

\subsection*{Reduced Size}
Our alternative to Chung's pebble extraction property was inspired by a detail in her proof. According to her proof, it is possible to extract at least $\lceil(\size{c}-\suppsize{c})/2\rceil$ pebbles from any configuration $c$ using pebbling steps of weight 2. How do we prove this, and how does it generalize to other weights? Extracting $n$ pebbles at a weight of $k$ is similar to finding $n$ blocks (pebbles placed on the same vertex) of $k$ pebbles. In the following theorem, we prove a lower bound for finding these blocks.

\begin{definition}
Let $V$ be a finite set.
Let $c:V\to\Nat$ and $n,k\in\Nat$. We say that $n$ blocks of $k$ can be extracted from $c$ if there is a function $e:V\to\Nat$ indicating the number of blocks that can be extracted from each vertex, such that $\size{e}=n$ and $c(v)\ge k\cdot e(v)$ for all $v\in V$. The result of this extraction is a configuration $c'$ defined as $c'(v)\defeq c(v)-k\cdot e(v)$.
\end{definition}

\begin{theorem}
\label{thm:block_ext}
Let $V$ be a finite set.
Let $c:V\to\Nat$ and $k\ge 1$.\\
If $\size{c}-(k-1)(\suppsize{c}-1)\ge nk$, then $n$ blocks of $k$ can be extracted from $c$.
\end{theorem}
\begin{proof}
Use induction on $n$. The proof is trivial for $n=0$. Suppose the theorem holds for~$n$, and $\size{c}-(k-1)(\suppsize{c}-1)\ge (n+1)\cdot k$. Note that $\size{c}>\suppsize{c}\cdot(k-1)$. Using the pigeonhole principle determine $v\in V$ with $c(v)\ge k$. Subtract $k$ from $c(v)$ and call the resulting configuration $c^\star$. Now $\size{c^\star}-(k-1)(\suppsize{c^\star}-1)\ge nk$, and using the induction hypothesis we can extract $n$ blocks of $k$ from $c^\star$. In total $n+1$ blocks of $k$ are extracted from~$c$, proving the theorem.
\end{proof}

\begin{corollary}
Let $V$ be a finite set.
Let $c:V\to\Nat$ and $k\in\Nat$.\\
If $n=\lceil(\size{c}-\suppsize{c}\cdot(k-1))/k\rceil$ then $n$ blocks of $k$ can be extracted from $c$.
\end{corollary}
\begin{proof}
Apply the previous theorem. Note that $\lceil x/k\rceil=\lfloor (x+k-1)/k\rfloor$ for any $x$.
\end{proof}

In Theorem \ref{thm:block_ext} you can think of $(k-1)(\suppsize{c}-1)$ as a kind of offset that accounts for vertices not having an exact multiple of $k$ pebbles. Any combination of blocks of at most $k$ pebbles that can be fitted into the configuration size that remains after subtracting this offset can be extracted, no matter how spread out the configuration is. If all pebbles are placed on a single vertex, such that $\suppsize{c}=1$, then no offset is needed. To use this concept, we define the $k$-reduced size of a configuration:

\begin{definition}
The $k$-reduced size $\mathsf{r}_k$ of a configuration $c$ is defined as:
$$\rsize{k}{c}\defeq\size{c}-(k-1)(\suppsize{c}-1)$$
\end{definition}

\begin{lemma}
Let $V$ be a finite set.
For $c_1,c_2:V\to\Nat$ holds
$\rsize{k}{c_1+c_2}\ge\rsize{k}{c_1}+\rsize{k}{c_2}$.
\end{lemma}
\begin{proof}
Left to the reader. Note that the support of $c_1$ and $c_2$ may intersect.
\end{proof}

\begin{lemma}
\label{lem:rsize_extract}
Let $V$ be a finite set.
Let $c:V\to\Nat$ and $k,n\in\Nat$.
If~\,$\rsize{k}{c}\ge n$ and $n\le k$,\linebreak then one block of $n$ pebbles can be extracted from $c$, resulting in $c'$ with $\rsize{k}{c'}\ge\rsize{k}{c}-n$.
\end{lemma}
\begin{proof}
Left to the reader. Note that the support may decrease.
\end{proof}

\section{The Induction Step}
Our alternative to property (ii) of Theorem \ref{thm:hypercube_2} uses the concept of reduced configuration size. Instead of directly extracting pebbles, we prove that two pebbles can be moved to a target vertex in such a way that the configuration which is left over has a certain $k$-reduced size. This will ensure that enough pebbles can be extracted from it later. We define this property in the form of a new kind of pebbling number written as $\tau_{n,k}$.

\begin{definition}
Let $G$ be a graph. Given configurations $c_1$ and $c_2$ on $G$, we say that $c_1$ is a subconfiguration of $c_2$, denoted as $c_1\subseteq c_2$, if $c_1(v)\le c_2(v)$ for all $v\in V(G)$.
\end{definition}

\begin{definition}
Let $G$ be a graph. The number $\taupebb{n}{k}{G,t}$ is the smallest $p\in\Nat$ where for every configuration $c$ with $\size{c}=p-\suppsize{c}+1+m$, there is a subconfiguration $c^\star$ of $c$ such that $c^\star$ is $n$-fold $t$-solvable, and the $k$-reduced size of the residual configuration $c-c^\star$ is at least\,\,$m$:
\begin{align*}
\taupebb{n}{k}{G,t}\defeq\min\Set{~p\in\Nat~}{~\begin{aligned}
&\forall m\in\Nat~\forall c:\size{c}=p-\suppsize{c}+1+m~\implies\\
&\exists c^\star\subseteq c:\solv{n}{t}{c^\star}~\land~\rsize{k}{c-c^\star}\ge m
\end{aligned}~}
\end{align*}
\end{definition}
\vspace{1mm}

\newcommand{\squareconf}[4]{\begin{bmatrix}#1&#2\\#3&#4\end{bmatrix}}
\paragraph{Example}
Let's illustrate the definition of $\tau_{n,k}$ on the graph $G=\pathgraph{2}{3}\gprod\pathgraph{2}{4}$. We will denote configurations on $G$ using a square matrix, where steps in the horizontal direction have weight~4 and steps in the vertical direction have weight 3. The target vertex $t$ will be the lower left corner. Later in this chapter it will be easy to see that the following holds:
$$\taupebb{2}{3}{\pathgraph{2}{3}\gprod\pathgraph{2}{4},~t}\le 2\cdot 3\cdot 4=24$$
According to this bound, every configuration $c$ on $G$ with $\size{c}=24-\suppsize{c}+1+m$ has a 2-fold $t$-solvable subconfiguration such that the residual configuration has a 3-reduced size of at least~$m$. If we are given a configuration that puts 39 pebbles in the top right corner, so that $m=15$, then we must pick the subconfiguration with 24 pebbles on this corner. This subconfiguration is 2-fold $t$-solvable, and the residual configuration has a 3-reduced size of~15 (equal to its normal size because the support is one), which is sufficient.
$$\rsize{3}{\squareconf{0}{39}{0}{0}-\squareconf{0}{24}{0}{0}}=15\ge m=15$$
If we are instead given a configuration that is more spread out, say 13 pebbles on each of the corners except $t$, then we have more freedom. Note that now $m=17$, since the support of this configuration contains 3 vertices. We could for example use one of the corners to form a subconfiguration that is 2-fold $t$-solvable, and still have plenty of pebbles left to make sure that the residual configuration has a sufficient 3-reduced size.
$$\rsize{3}{\squareconf{13}{13}{0}{13}-\squareconf{13}{0}{0}{0}}=24\ge m=17$$

\pagebreak
\begin{lemma}
\label{lem:tau_K1}
If $v$ is the only vertex of $K_1$, then:
$$\taupebb{2}{k}{K_1,v}\le 2$$
\end{lemma}
\begin{proof}
Let $c$ be a configuration on $K_1$ with $\size{c}=2-\suppsize{c}+1+m=2+m$. Define the subconfiguration $c^\star$ as $c^\star(v)=2$. Now $c^\star$~is 2-fold $v$-solvable and $\rsize{k}{c-c^\star}=m$.
\end{proof}

\begin{lemma}
\label{lem:taupebb_extend}
Let $G$ be a graph.
If $\pebb{G,t}\le p$ and $\taupebb{2}{k}{G,t}\le 2p$, then for $n\ge 2$:
$$\taupebb{n}{k}{G,t}\le np$$
\end{lemma}
\begin{proof}
Let $c$ be a configuration on $G$ with $\size{c}=np-\suppsize{c}+1+m$. First move two pebbles to $t$ using $\taupebb{2}{k}{G,t}\le 2p$, and then move another $n-2$ pebbles to $t$ from the residual configuration using $\pebb{G,t}\le p$. The details are left to the reader.
\end{proof}

Our goal is to prove Theorem \ref{thm:hypercube_2} for hypercubes with generalized edge weights. We prove the induction step as two separate lemmas. Instead of two $n$-cubes forming an $(n+1)$-cube, we will have two copies of an arbitrary graph with edges between matching vertices. To record the fact that between these two identical graphs, which were previously called $M_1$ and $M_2$, pebbles are only moved in one direction, we define the \emph{arrow graph}; a single directed edge with weight $k$ going from 0 to 1.

\begin{definition}
Let $k\in\Nat$. The arrow graph denoted by $\arrowgraph{k}$ is defined as $(V,E,\weight)$ where: $$V=\{\arrowfrom,\arrowto\},\quad E=\{(\arrowfrom,\arrowto)\},\quad\weight(\arrowfrom,\arrowto)=k.$$
\end{definition}

\SimpleFigure{fig/example/Cartesian_product_arrows.svg}{The Cartesian product of 3 arrow graphs.}

\pagebreak
\SimpleFigure{fig/example/matching_graphs.svg}{}

\begin{lemma}
\label{lem:pi_arrow_prod}
Let $G$ be a graph and $t$ a target vertex.
Let $k,l\in\Nat$ such that $2\le k\le l$.\\
If $\pebb{G,t}\le p$ and $\taupebb{2}{l}{G,t}\le 2p$ then:
$$\pebb{\arrowgraph{k}\gprod G,~(\arrowto,t)}\le kp$$
\end{lemma}
\begin{proof}
Let $c$ be a configuration on $\arrowgraph{k}\gprod G$ such that $|c|\ge kp$. Split the graph into $G_0$ containing the vertex $(\arrowfrom,t)$ and $G_1$ containing the vertex $(\arrowto,t)$. Split $c$ into the configurations $c_0$ on $G_0$ and $c_1$ on $G_1$. Consider the following two cases:

\vspace{2mm}
\item\textbf{(i)}~~$\size{c_1}\ge\suppsize{c_0}$.
\begin{itemize}[--,leftmargin=*]
\item Use Lemma \ref{lem:rsize_extract} to move $\lfloor\rsize{k}{c_0}/k\rfloor$ pebbles to $G_1$. Note that $k\le l$.
\item Derive $\size{c_1}+\lfloor\rsize{k}{c_0}/k\rfloor\ge p$ from the following inequalities:
\begin{align*}
\rsize{k}{c_0}&\le k\lfloor\rsize{k}{c_0}/k\rfloor+(k-1)\\
k\size{c_1}+\rsize{k}{c_0}&=k\size{c_1}+\size{c_0}-(k-1)(\suppsize{c_0}-1)\\
&\ge k\size{c_1}-(k-1)\size{c_1}+(k-1)+\size{c_0}\\
&\ge kp+(k-1)
\end{align*}
\item Now $G_1$ contains at least $p$ pebbles, which is enough to move one pebble to $(\arrowto,t)$.
\end{itemize}
\vspace{2mm}
\item\textbf{(ii)}~~$\size{c_1}\le\suppsize{c_0}-1$.
\begin{itemize}[--,leftmargin=*]
\item Derive $\size{c_0}\ge kp-\size{c_1}\ge kp-\suppsize{c_0}+1$.
\item Use Lemma \ref{lem:taupebb_extend} to move $k$ pebbles to $(\arrowfrom,t)$, and from there move one pebble to $(\arrowto,t)$.
\end{itemize}
\end{proof}

\pagebreak
\begin{lemma}
\label{lem:tau_arrow_prod}
Let $G$ be a graph and $t$ a target vertex.
Let $k,l\in\Nat$ such that $2\le k\le l$.\\
If $\pebb{G,t}\le p$ and $\taupebb{2}{l}{G,t}\le 2p$ then:
$$\taupebb{2}{k}{\arrowgraph{k}\gprod G,~(\arrowto,t)}\le 2kp$$
\end{lemma}
\begin{proof}
Let $c$ be a configuration on $\arrowgraph{k}\gprod G$ such that $|c|=2kp-\suppsize{c}+1+m$. Again\linebreak split the graph into $G_0$ and $G_1$, and split $c$ into $c_0$ and $c_1$. Consider the following two cases:

\item\textbf{(i)}~~$\size{c_1}\le 2p-\suppsize{c_1}$.
\begin{itemize}[--,leftmargin=*]
\item Define $m^\star\defeq kp-\size{c_1}-\suppsize{c_1}$.
\item Derive $\size{c_0}\ge kp-\suppsize{c_0}+1$:
\begin{align*}
\size{c_0}&\ge 2kp-\suppsize{c}+1-\size{c_1}\\
&=(kp-\suppsize{c_0}+1)+(kp-\suppsize{c_1})-\size{c_1}\\
&\ge kp-\suppsize{c_0}+1
\end{align*}
\item Move $k$ pebbles to $(\arrowfrom,t)$ in $G_0$ using Lemma \ref{lem:taupebb_extend}, resulting in a residual configuration $c_0^\ast$ with $\rsize{k}{c_0^\ast}\ge m^\star+m$:
\begin{align*}
\rsize{k}{c_0^\ast}&\ge{\size{c_0}-kp+\suppsize{c_0}-1}\\
&=kp+\size{c_0}+\suppsize{c_0}-2kp-1\\
&=kp+\size{c}-\size{c_1}+\suppsize{c}-\suppsize{c_1}-2kp-1\\
&=m^\star+m
\end{align*}
\item Move $\lfloor m^\star/k\rfloor$ pebbles to $G_1$ using $c_0^\ast$ and derive $\size{c_1}+\lfloor m^\star/k\rfloor\ge p$. What remains of $c_0^\ast$ has a $k$-reduced size of at least $m$.
\item There are now at least $p$ pebbles on $G_1$, so one pebble can be moved to $(\arrowto,t)$ in $G_1$, and another pebble can be moved to $(\arrowto,t)$ using the $k$ pebbles that were moved to $(\arrowfrom,t)$.
\end{itemize}

\item\textbf{(ii)}~~$\size{c_1}\ge 2p-\suppsize{c_1}+1$.
\begin{itemize}[--,leftmargin=*]
\item Move two pebbles to $(\arrowto,t)$ in $G_1$ using the hypothesis $\taupebb{2}{l}{G,t}\le 2p$, resulting in a residual configuration $c_1^\ast$.
\item Derive $\rsize{k}{c_0+c_1^\ast}\ge m$:
\begin{align*}
\suppsize{c_0}&\le\Count{V(G)}\le p\\
\suppsize{c_0}&\le(k-1)(2p-\suppsize{c_0}+1)-1\\
\rsize{k}{c_0+c_1^\ast}&\ge\rsize{k}{c_0}+\rsize{k}{c_1^\ast}\\
&\ge\size{c_0}-(k-1)(\suppsize{c_0}-1)+\size{c_1}-2p+\suppsize{c_1}-1\\
&=\size{c}-2kp+\suppsize{c_1}+(k-1)(2p-\suppsize{c_0}+1)-1\\
&\ge\size{c}-2kp+\suppsize{c}=m
\end{align*}
\end{itemize}
\end{proof}

\pagebreak
\section{Results}
\begin{theorem}
For~~$2\le k_1\le k_2\le\dots\le k_n$:
$$\pebb{\arrowgraph{k_1}\gprod\arrowgraph{k_2}\gprod\dots\gprod\arrowgraph{k_n},(\arrowto,\arrowto,\dots,\arrowto)}={k_1}{k_2}\dots{k_n}$$
\end{theorem}
\begin{proof}
This theorem results from inductively applying Lemma \ref{lem:pi_arrow_prod} and \ref{lem:tau_arrow_prod}. To prove the base case, note that $\arrowgraph{k_n}$ is isomorphic to $\arrowgraph{k_n}\gprod K_1$, and use Lemma \ref{lem:tau_K1}.
To see that the bound obtained so far is strict, put $k_1 k_2\dots k_n-1$ pebbles on $(\arrowfrom,\arrowfrom,\dots,\arrowfrom)$ and try to move one pebble to $(\arrowto,\arrowto,\dots,\arrowto)$.
\end{proof}

\begin{corollary}
\label{thm:hypercube}
For~~$2\le k_1\le k_2\le\dots\le k_n$:
$$\pebb{\pathgraph{2}{k_1}\gprod\pathgraph{2}{k_2}\gprod\dots\gprod\pathgraph{2}{k_n}}= {k_1}{k_2}\dots{k_n}$$
\end{corollary}
\begin{proof}
Use the previous theorem.
\end{proof}

\subsection*{Generalization}
Chung says that $\pebb{K_n\gprod G}=n\cdot\pebb{G}$ when all edge weights are 2 and $G$ has the\linebreak 2-pebbling property. We can state a similar generalization of the results we just presented using a star-like digraph with all edges pointing towards the central vertex. We denote this graph as $\instargraph{n}{k}$, where $n$ is the number of leafs and every edge has weight $k$. Note that the arrow graph $\arrowgraph{k}$ is equivalent to $\instargraph{1}{k}$. We label the central vertex of $\instargraph{n}{k}$ as $s$. The pebbling number of $s$ can be determined using a proof similar to Theorem \ref{thm:pebb_K}:
$$\pebb{\instargraph{n}{k},~s}=nk-n+1$$
Is the following also true?
\begin{conjecture}
Let $G$ be a graph and $n\ge 1$.
Let $k,l\in\Nat$ such that $2\le k\le l$.\\
If $\pebb{G,t}\le p$ and $\taupebb{2}{l}{G,t}\le 2p$ then:
\begin{enumerate}[(i)]
\item $\pebb{\instargraph{n}{k}\gprod G,~(s,t)}\le (nk-n+1)\cdot p$
\item $\taupebb{2}{k}{\instargraph{n}{k}\gprod G,~(s,t)}\le (nk-n+1)\cdot 2p$
\end{enumerate}
\end{conjecture}

\SimpleFigure{fig/graph/inward_star.svg}{Inward pointed star graphs}

\pagebreak

\chapter{Zero-sums}
\label{chapter:zero-sums}

With the results that we obtained in the previous chapter, we are finally ready to prove the \Erdos-Lemke conjecture (Theorem \ref{thm:erdos_lemke_conj}). To do this, we first show that pebbling solutions on one graph can be translated to solutions on another graph via a homomorphism. We use this approach to determine an upper bound for the pebbling number of a divisor lattice. Then we associate pebbles with subsequences, and use the pebbling steps to carry out a construction with these subsequences, which will finally produce the zero-sum subsequence we are looking for. Before starting, we note the following:

\begin{theorem}
\label{thm:php_zero_sum}
Every non-empty sequence of integers $a_1,a_2,\dots,a_n$ has a non-empty sub-sequence $b_1,b_2,\dots,b_m$ such that:
$$b_1+b_2+\dots+b_m\equiv 0~(\mathrm{mod}~n)$$
\end{theorem}
\begin{proof}
Define $c_0\defeq 0$ and $c_i\defeq(a_1+a_2+\dots+a_i)~\mathrm{mod}~n$.
Use the pigeonhole principle to determine $i$ and $j$ with $0\le i<j\le n$ such that $c_i=c_j$.
The subsequence $a_{i+1},a_{i+2},\dots,a_j$ has the desired property.
\end{proof}

\section{Homomorphisms}
A homomorphism from a graph $G$ to a graph $H$ is a mapping $\phi$ from vertices of $G$ to vertices of $H$ that is \emph{structure-preserving}, meaning that for every edge $(u,v)$ in $G$ there is also an edge $(\phi(u),\phi(v))$ in $H$. It is not necessary that $H$ contains a copy of $G$; sometimes it is possible to map multiple vertices in $G$ to one vertex in $H$. For example, for all bipartite graphs there is a structure-preserving mapping to $K_2$. We also want the homomorphism to preserve the edge weight, leading us to the following definition:

\begin{definition}
\label{def:homomorphisms}
Let $G$ and $H$ be graphs. A function $\phi:V(G)\to V(H)$ is a homomorphism from $G$ to $H$, also written as $\Hom{G}{H}{\phi}$, if for all $(u,v)\in E(G)$ we have:
$$(\phi(u),\phi(v))\in E(H)\qquad\text{and}\qquad\weight_G(u,v)=\weight_H(\phi(u),\phi(v)).$$
\end{definition}

When a homomorphism $\Hom{G}{H}{\phi}$ is surjective, we can use it to map a configuration $c$ on $H$ to a (not necessarily unique) configuration $c^\star$ on $G$ in the following way: for every $v\in V(H)$, determine a vertex $u\in V(G)$ such that $\phi(u)=v$, and set $c^\star(u)=c(v)$. Now a solution of  $c^\star$ on $G$ can be translated into a solution of $c$ on $H$ via the homomorphism.

\pagebreak
\begin{lemma}
\label{lem:pebb_hom}
Let $G$ and $H$ be graphs. If a surjective homomorphism $\Hom{G}{H}{\phi}$ exists, then:
\begin{enumerate}[(i)]
\item $\pebb{H}\le\pebb{G}$.
\item $\pebb{H,\phi(t)}\le\pebb{G,t}$ for all $t\in V(G)$.
\item $\taupebb{n}{k}{H,\phi(t)}\le\taupebb{n}{k}{G,t}$ for all $t\in V(G)$ and $n,k\in\Nat$.
\end{enumerate}
\end{lemma}
\begin{proof}
Left to the reader. When proving (iii) it should be noted that the $k$-reduced size of a configuration on $G$ does not decrease in the image of $\phi$.
\end{proof}

\subsection*{Grids}
In this section the vertices of a path $\pathgraph{n}{k}$ are treated as numbers from $0$ to $n-1$. We denote the $n$-dimensional hypercube as $Q_n$ and label its vertices with binary strings of zeros and ones. We define a special notation to count the number of ones in a binary string:
\begin{definition}
Let $x$ be a binary string. The number of ones in $x$ is denoted by $\bincount{x}$.
\end{definition}
\begin{theorem}
For~~$2\le k_1\le k_2\le\dots\le k_s$ and $n_1,n_2,\dots,n_s$ with $n_i\ge 1$ for all $i\le s$:
$$\pebb{\pathgraph{n_1}{k_1}\gprod\pathgraph{n_2}{k_2}\gprod\dots\gprod\pathgraph{n_s}{k_s}}= {k_1}^{n_1-1}{k_2}^{n_2-1}\dots{k_s}^{n_s-1}$$
\end{theorem}
\begin{proof}
Define a homomorphism $\phi$ as follows:
\begin{align*}
&\phi:Q_{n_1-1}^{(k_1)}\gprod Q_{n_2-1}^{(k_2)} \gprod\dots\gprod Q_{n_s-1}^{(k_s)}\longrightarrow\pathgraph{n_1}{k_1}\gprod\pathgraph{n_2}{k_2}\gprod\dots\gprod\pathgraph{n_s}{k_s}\\
&\phi((x_1,x_2,\dots,x_s))\defeq(\bincount{x_1},~\bincount{x_2},~\dots,~\bincount{x_s})
\end{align*}
Note that $\phi$ is a valid homomorphism, preserving both edge connectivity and weight, and that it is surjective. Using Lemma \ref{lem:pebb_hom} and Corollary \ref{thm:hypercube} we obtain an upper bound for the pebbling number of the product of paths. It is left to the reader to see that it is strict.
\end{proof}

\begin{figure}[h]
\includesvg{fig/example/hypercube_retraction.svg}
\vspace{3mm}
\caption{Retracting a hypercube to a $2\times 2$ grid.}
\end{figure}

\pagebreak
\subsection*{Divisor Lattices}
The divisor lattice of a number $n$ consists of all divisors of $n$ as vertices, including $n$ itself, and for each divisor an edge to every other divisor that is a multiple of it (excluding self-loops). In our definition, the edge weights are determined by the multiplication factor.

\begin{definition}
Let $n\in\Nat$. The divisor lattice $D_n$ is a graph $(V,E,\weight)$ where:
\begin{align*}
V&=\Set{d\in\Nat}{(d\mid n)}\\
E&=\Set{(p,d)\in V\times V}{p\neq d\land (p\mid d)}\\
\weight(p,d)&= d/p
\end{align*}
\end{definition}

We want to know the pebbling number of this graph. By looking at the prime factorization of $n$'s divisors, we can see that $D_n$ contains a grid with the length of each dimension being the exponent of one of $n$'s prime factors plus one. Using the result from the previous section we can then determine that $\pebb{D_n,n}\le n$. This is essentially the approach taken by Chung. To make it easier to formalize, we instead prove this result directly from the step lemmas that we also used in the previous section.

\begin{theorem}
Let $n\ge 2$. If $p$ is the smallest prime factor of $n$, then:
$$\pebb{D_n,n}\le n\qquad\text{and}\qquad\taupebb{2}{k}{D_n,n}\le 2n.$$
\end{theorem}
\begin{proof}
We prove this theorem by induction on $n$. Suppose the theorem has been proven for all $m<n$. Let $m\defeq n/p$, and define a homomorphism $\phi:\arrowgraph{p}\gprod D_m\to D_n$ as follows:
\begin{align*}
\phi((\arrowfrom,d))&\defeq d\\
\phi((\arrowto,d))&\defeq pd
\end{align*}
Note that $\phi$ is well-defined and surjective. Using Lemma \ref{lem:pebb_hom} we find:
\begin{align*}
\pebb{D_n,n}&\le\pebb{\arrowgraph{p}\gprod D_m,~(\arrowto,m)}\\ \taupebb{2}{k}{D_n,n}&\le\taupebb{2}{k}{\arrowgraph{p}\gprod D_m,~(\arrowto,m)}
\end{align*}
Now use Lemma \ref{lem:pi_arrow_prod} and \ref{lem:tau_arrow_prod}, the induction steps for $\pi$ and $\tau_{2,k}$ from the previous chapter. If $m\ge 2$ then the proof can be finished using the induction hypothesis. Note that $p\le p_m$ where $p_m$ is the smallest prime factor of $m$.

Otherwise, if $m=1$, use that $D_1$ is isomorphic to $K_1$.
\end{proof}

\pagebreak
\section{Pebbling Construction}
Now we can finally show how graph pebbling is connected to number theory. Pebbles will carry sequences of numbers, and a predicate will determine if a sequence is \emph{well-placed} on a certain vertex, such that a pebble carrying that sequence may be put there. A pebbling step will represent a way to pick a subsequence from the union of a number of sequences. Moving one pebble to a target vertex $t$ will result in the construction of a subsequence with the property that it is well-placed on $t$.

\newcommand{\wellplaced}[1]{\mathcal{P}_{#1}}
\begin{lemma}
\label{lem:pebbling_construction}
Let $G$ be a graph and $A$ an set. Let $a_1,a_2,\dots,a_n\in A$ be a sequence of elements in $A$. Define $I\defeq\{1,2,\dots,n\}$. For every $v\in V(G)$, let a predicate $\wellplaced{v}$ on subsets of $I$ be given that defines if a subsequence $\{a_i\}_{i\in S}$, determined by a set $S\subseteq I$, is well-placed on $v$. Suppose:
\begin{enumerate}[(I)]
\item For every $i\in I$ there is a vertex $v$ such that $\wellplaced{v}(\{i\})$.
\item For every edge $(u,v)$ in $G$ and every sequence $S_1,S_2,\dots,S_k$ of disjoint and non-empty subsets of $I$, where $k=\weight_G(u,v)$ and $\wellplaced{u}(S_i)$ for all $i\le k$, there exists a non-empty\linebreak set $S\subseteq S_1\cup S_2\cup\dots\cup S_k$ such that $\wellplaced{v}(S)$.
\end{enumerate}
\noindent If $\pebb{G,t}\le n$ for a vertex $t$, then there exists a non-empty set $S\subseteq I$ such that $\wellplaced{t}(S)$.
\end{lemma}
\begin{proof}
Put every $a_i$ on a vertex determined by (I). Define an initial pebbling configuration $c$ by counting how many $a_i$'s are placed on each vertex, and note that $c$ is $t$-solvable since $\size{c}=n$. Using (II) and a solution of~$c$ that puts one pebble on $t$, we construct a subsequence of $a_1,a_2,\dots,a_n$ that is well-placed on $t$.
\end{proof}

\subsection*{Divisor Pebbles}
To demonstrate the concept of a pebbling construction more clearly, we are first going to prove the \Erdos-Lemke conjecture for $n=d$:

\begin{theorem}
Let $n\ge 1$. Every sequence $a_1,a_2,\dots,a_n$ of $n$ divisors of $n$ has a non-empty subsequence $b_1,b_2,\dots,b_m$ such that:
$$b_1+b_2+\dots+b_m=n$$
\end{theorem}
\begin{proof}
Use Lemma \ref{lem:pebbling_construction} with $G=D_n$ and the sequence $a_1,a_2,\dots,a_n$. Let a subsequence be well-placed on a vertex $d\in V(D_n)$ if its sum is equal to $d$. For each $i\le n$ we can place the one-element sequence $a_i$ on the vertex $a_i$. If a vertex $d\in V(D_n)$ has at least $p$ pebbles, in the form of $p$ disjoint subsequences that each sum to $d$, then the union of these subsequences sums to $pd$, allowing us to put one pebble on the vertex $pd$ in exchange for removing $p$ pebbles from vertex $d$. This proves (I) and (II) of Lemma \ref{lem:pebbling_construction}. Since $\pebb{D_n,n}\le n$, it follows that there is a subsequence of $a_1,a_2,\dots,a_n$ which sums to~$n$.
\end{proof}

\pagebreak
\paragraph{Example} We consider the case $n=60=2^2\cdot3\cdot5$. The configuration shown below, on the Hasse diagram of $D_{60}$, is one possible way of picking 60 divisors of 60: 29 ones, 15 twos, 10~threes, and 6 fives. The edge weight is 2 for going left, 3 for going up, and 5 for going right. A solution of this configuration, for target vertex 60, is shown as a pebble flow. We put 5 pebbles on 6 using 15 twos, and another 5 using 10 threes, then we use the 10 pebbles on 6 to put 2 pebbles on 30, and finally we can put one pebble on 60.

\vspace{2cm}
\begin{figure}[h]
\includesvg{fig/example/flow_divisors.svg}
\vspace{1cm}
\caption{}
\end{figure}

\newpage
\subsection*{GCD Pebbles}
To prove the case $n>d$, Lemke and Kleitman used a generalization of the theorem based on \emph{greatest common divisors}. We will prove this generalization following the approach of Chung, including a technical correction described in \cite{Clarke1997}. We made a computer formalization of this proof using the \emph{Coq Proof Assistant}, following the same steps as described in this and the previous chapter. The formalization required investigating all reasoning steps in detail, and gives an almost absolute guarantee that the proof, which we are about to finish, is correct. More details about the formalization can be found in the appendix.

\begin{theorem}
Let $n\ge 1$. For every sequence of positive integers $a_1,a_2,\dots,a_n$ there is a non-empty subset $S$ of $\{1,2,\dots,n\}$ such that:
$$\sum_{i\in S}a_i\equiv 0~(\mathrm{mod}~n)\qquad\text{and}\qquad\sum_{i\in S}\gcd(n,a_i)\le n.$$
\end{theorem}
\begin{proof}
We again use Lemma \ref{lem:pebbling_construction}. This time, a subsequence determined by $S\subseteq \{1,2,\dots,n\}$ is well-placed at vertex $d\in V(D_n)$ if:
\begin{enumerate}[(a)]
\item $\sum_{i\in S}a_i$ is a multiple of $d$,
\item $\sum_{i\in S}\gcd(n,a_i)\le d$.
\end{enumerate}
It is now sufficient to prove property (I) and (II) of Lemma \ref{lem:pebbling_construction}:
\begin{enumerate}[(I)]
\item For all $i\le n$ the one-element sequence $a_i$ is well-placed at vertex $\gcd(n,a_i)$.
\item Let $(d,pd)$ be an edge of $D_n$. Note that $\weight(d,pd)=p$. Suppose a sequence of $p$ disjoint sets $S_1,S_2,\dots,S_p\subseteq\{1,2,\dots,n\}$ is given that are all well-placed at vertex $d$. Define:
\begin{align*}
x_i&\defeq\textstyle\sum_{j\in S_i}a_j\\
r_i&\defeq x_i/d
\end{align*}
Note that $x_i$ is a multiple of $d$ according to (a). Using Theorem \ref{thm:php_zero_sum}, determine a set $R\subseteq\{1,2,\dots,p\}$ corresponding to a subsequence of $r_1,r_2,\dots,r_p$ that sums to a multiple of $p$. Now the set $S$ given by $\bigcup_{i\in R}S_i$ is well-placed at vertex $pd$.
\end{enumerate}
\end{proof}

\pagebreak
\begin{corollary}
\emph{(\Erdos-Lemke conjecture)} Given positive integers $n$, $d$, and $a_1,a_2,\dots,a_d$ with $d\,|\,n$ and $a_i\,|\,n$ for all $i\le d$, there is a non-empty subset $S$ of $\{1,2,\dots,d\}$ such that:
$$\sum_{i\in S}a_i\equiv 0~(\mathrm{mod}~d)\qquad\text{and}\qquad\sum_{i\in S}a_i\le n.$$
\end{corollary}
\begin{proof}
Apply the previous theorem with $n\defeq d$. For the second requirement, use:
$$\sum_{i\in S}a_i=\sum_{i\in S}\gcd(n,a_i)\le\sum_{i\in S}{n/d}\cdot\gcd(d,a_i)\le n$$
This is based on the following property of $\gcd$ with $k=n/d$:
$$\gcd(k\cdot d,a)\le k\cdot\gcd(d,a)$$
\end{proof}

\vfill
\begin{center}
\includegraphics[width=.3\textwidth]{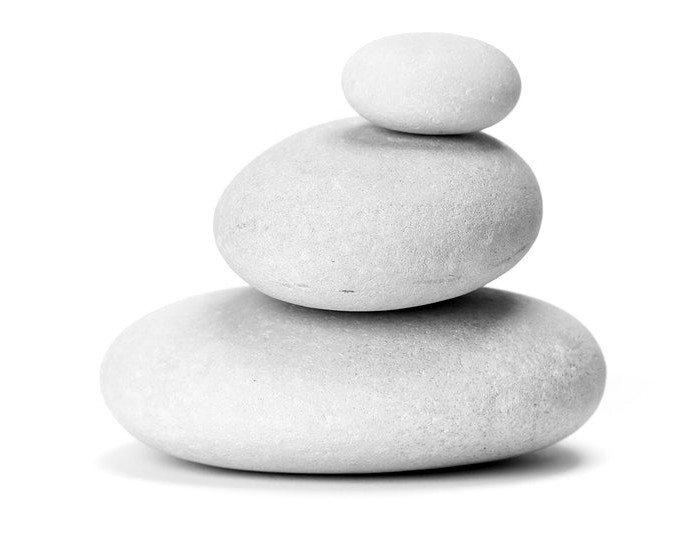}
\end{center}

\newpage
\begin{figure}
\includegraphics[width=.72\textwidth]{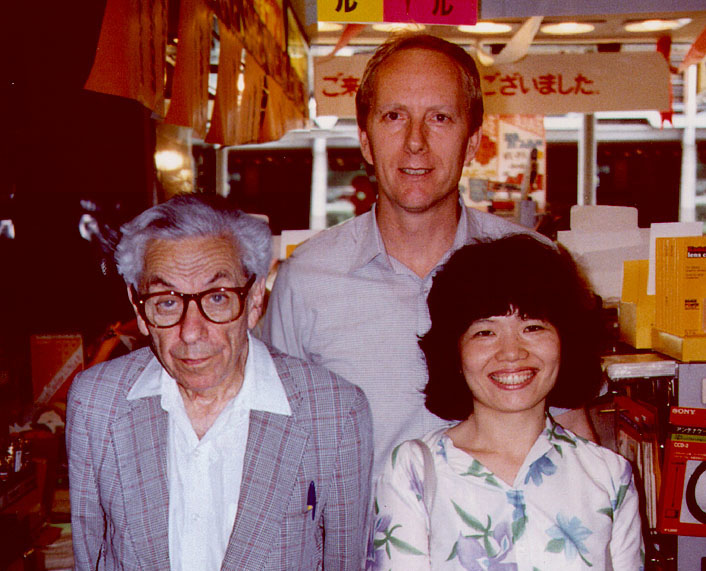}
\vspace{2mm}
\caption{Paul \Erdos, Ronald Graham and Fan Chung (1986).}
\end{figure}

\cleardoublepage
\phantomsection
\addcontentsline{toc}{chapter}{Bibliography}
\setlength\bibitemsep{3\itemsep}
\printbibliography

%
%
%

\makeatletter
\newcommand{\replunderscores}[1]{\expandafter\@repl@underscores#1_\relax}
\def\@repl@underscores#1_#2\relax{%
    \ifx \relax #2\relax
        #1%
    \else
        #1%
        \_
        \@repl@underscores#2\relax
    \fi
}
\makeatother

\newcommand{\tok}[1]{\mathrm{\replunderscores{#1}}}
\newcommand{\Coq}[1]{\mathbf{#1}}
\newcommand{\CoqVariable}{\Coq{Variable}}
\newcommand{\CoqDefinition}{\Coq{Definition}}
\newcommand{\CoqHypothesis}{\Coq{Hypothesis}}
\newcommand{\CoqLemma}{\Coq{Lemma}}
\newcommand{\CoqTheorem}{\Coq{Theorem}}
\newcommand{\CoqType}{\tok{Type}}
\newcommand{\CoqProp}{\tok{Prop}}
\newcommand{\coqnat}{\tok{nat}}
\newcommand{\coqlist}{\tok{list}}
\newcommand{\coqlistnat}{\tok{list}~\tok{nat}}

\renewcommand{\implies}{~\rightarrow~}
\newcommand{\submseteq}{\subseteq^+}
\newcommand{\fmap}{\mathbin{\,\raisebox{1pt}{\scalebox{0.8}{$<$}}\$\raisebox{1pt}{\scalebox{0.8}{$>$}}\,}}

\chapter*{Appendix: Coq Formalization}
\addcontentsline{toc}{chapter}{Appendix}
\markboth{Coq Formalization}{Coq Formalization}

The purpose of this appendix is to show some of the highlights of our Coq\footnote{\url{https://coq.inria.fr/}} formalization of the \Erdos-Lemke conjecture. The source code of this formalization is available online\footnote{\url{https://github.com/bergwerf/pebbling}, \textsc{doi:} \href{https://doi.org/10.5281/zenodo.7876086}{10.5281/zenodo.7876086}}. Here we use some slightly different notations than in the Coq source to improve readability.
The formalization was developed using version \emph{8.16.0} of Coq, and version \emph{1.8.0} of the the \emph{Coq-std++} library\footnote{\url{https://gitlab.mpi-sws.org/iris/stdpp}}. Future versions of these systems may introduce breaking changes, and so it is unsure if our formalization can still be evaluated a few decades from now.

\begin{figure}[h]
\begin{tcolorbox}[ams nodisplayskip]
\begin{flalign*}
&\Coq{Corollary}~\tok{Erdos_Lemke_conjecture}~(l:\coqlistnat)~(d~n:\coqnat):&&\\
&\quad\tok{length}~l=d\implies d>0\implies n>0\implies
d\mid n\implies(\forall a,~a\in l\implies a\mid n) \implies&&\\
&\quad\exists~(l':\coqlistnat),~l'\neq []\land l'\submseteq l~\land
d\mid\tok{summation}~l'~\land
\tok{summation}~l'\le n.&&
\end{flalign*}
\end{tcolorbox}
\caption{The \Erdos-Lemke conjecture in Coq.}
\end{figure}

\paragraph{Naming Heuristics}
Like in most programming languages, all definitions and theorems in Coq must have a unique identifier expressed as a plain text string with various restrictions. Finding a naming pattern for these identifiers that is both coherent and convenient is a significant challenge when developing a large computer formalization. One of the naming \emph{heuristics} (it is too loose to be called a convention) that we used to form new identifiers is appending the main identifiers in the new definition or theorem in their order of appearance. We proved many small lemmas about summations, of which three simple examples are shown below. Here $\submseteq$ denotes the multiset inclusion operator from Coq-std++.

\begin{figure}[h]
\begin{tcolorbox}[ams nodisplayskip]
\begin{flalign*}
&\CoqLemma~\tok{elem_of_summation}:n\in l\implies n\le\tok{summation}~l.&&\\
&\CoqLemma~\tok{submseteq_summation}:l_1\submseteq l_2\implies\tok{summation}~l_1\le\tok{summation}~l_2.&&\\
&\CoqLemma~\tok{summation_nonzero}:\tok{summation}~l>0\implies\exists~n,~n\in l\land n>0.&&
\end{flalign*}
\end{tcolorbox}
\caption{Lemmas about the summation function.}
\end{figure}

\paragraph{Pebbling Bounds}
In the formalization it is not convenient to define the pebbling number as a function. We instead use a predicate that defines when, for a given target vertex, the pebbling number is \emph{bounded} by a certain value:
\begin{itemize}
\item $\pebb{G,t}\le p$ is defined by: $\tok{vertex_pebbling_bound}~G~p~t$,
\item $\taupebb{n}{k}{G,t}\le p$ is defined by: $\tok{vertex_pebbling_property}~G~n~k~p~t$.
\end{itemize}
Using these predicates, we can formulate the main lemmas from Chapter \ref{chapter:hypercubes} in Coq as shown in the excerpt below. Here the graph called `arrow' is analogous to the arrow graph $\arrowgraph{k}$.

\begin{figure}[h]
\begin{tcolorbox}[ams nodisplayskip]
\begin{flalign*}
&\CoqHypothesis~H_1:2\le k\le l.&&\\
&\CoqHypothesis~H_2:\tok{vertex_pebbling_bound}~G~p~t.&&\\
&\CoqHypothesis~H_3:\tok{vertex_pebbling_property}~G~2~l~(2\cdot p)~t.&&\\[3mm]
&\CoqLemma~\tok{vertex_pebbling_bound_arrow_prod}:&&\\
&\quad\tok{vertex_pebbling_bound}~(\tok{arrow}\gprod G)~(k\cdot p)~(v_1,t)&&\\[3mm]
&\CoqLemma~\tok{vertex_pebbling_property_arrow_prod}:&&\\
&\quad\tok{vertex_pebbling_property}~(\tok{arrow}\gprod G)~2~k~(2\cdot k\cdot p)~(v_1,t)&&
\end{flalign*}
\end{tcolorbox}
\caption{Results analogous to Lemma \ref{lem:pi_arrow_prod} and \ref{lem:tau_arrow_prod}.}
\end{figure}

\paragraph{Homomorphisms}
Transferring pebbling bounds between graphs using surjective homomorphisms was surprisingly difficult to formalize. The main difficulty is the inverse mapping of configurations that we described after Definition \ref{def:homomorphisms}.

\begin{figure}[h]
\begin{tcolorbox}[ams nodisplayskip]
\begin{flalign*}
&\CoqHypothesis~\tok{surj}:\tok{Surj}~({=})~h.&&\\
&\CoqHypothesis~\tok{hom}:\tok{Graph_Hom}~G~H~h.&&\\[3mm]
&\CoqLemma~\tok{surj_hom_vertex_pebbling_bound}~n~t:&&\\
&\quad\tok{vertex_pebbling_bound}~G~n~t\implies\tok{vertex_pebbling_bound}~H~n~(h~t).&&
\end{flalign*}
\end{tcolorbox}
\caption{A result similar to Lemma \ref{lem:pebb_hom}(ii).}
\end{figure}

\pagebreak
\noindent In Coq, we defined a small functional program called `inv\_dmap' (\emph{inverse domain mapping}) that carries out this inversion procedure. Among other properties, we had to prove that the size and the number of elements in the support remain unchanged when this function is applied to a configuration.

\begin{figure}[h]
\begin{tcolorbox}[ams nodisplayskip]
\begin{flalign*}
&\CoqDefinition~\tok{inv}~(f:A\to B)~(b:B):\coqlist~A\defeq&&\\
&\quad\tok{filter}~(\lambda~a,~f~a=b)~(\tok{enum}~A).&&\\[3mm]
&\CoqDefinition~\tok{dmap}~(f:A\to\coqnat):B\to\coqnat\defeq&&\\
&\quad\lambda~b,~\tok{summation}~(f\fmap\tok{inv}~h~b).&&\\[3mm]
&\CoqDefinition~\tok{inv_dmap}~(f:B\to\coqnat):A\to\coqnat\defeq&&\\
&\quad\lambda~v,~\tok{if}~\tok{decide}~(\tok{head}~(\tok{inv}~h~(h~v))=\tok{Some}~v)~\tok{then}~f~(h~v)~\tok{else}~0.&&
\end{flalign*}
\end{tcolorbox}
\caption{Configuration mappings.}
\end{figure}

\paragraph{Pebbling Construction}
The construction described in Lemma \ref{lem:pebbling_construction} is carried out in Coq using lists of elements of an arbitrary type $A$. A list called `pebbles' determines the values that may be used. Each vertex in the graph is associated with any number of lists with elements of type $A$, and the concatenation of all these lists should be a permutation of (a part of) the pebbles list.

\begin{figure}[h]
\begin{tcolorbox}[ams nodisplayskip]
\begin{flalign*}
&\CoqVariable~\tok{pebbles}:\coqlist~A.&&\\
&\CoqVariable~\tok{well_placed}:V~G\to\coqlist~A\to\CoqProp.&&\\[3mm]
&\CoqHypothesis~\tok{initialize}:\forall~a,~a\in\tok{pebbles}\implies\exists~v,~\tok{well_placed}~v~[a].&&\\
&\CoqHypothesis~\tok{step}:\forall~(u~v:V~G)~(l^\ast:\coqlist~(\coqlist~A)),&&\\
&\quad E~G~u~v\implies\tok{Forall}~(\tok{well_placed}~u)~l^\ast\implies
\tok{length}~l^\ast=\tok{weight}~G~u~v\implies&&\\
&\quad \exists~(l:\coqlist~A),~l\submseteq\tok{concat}~l^\ast\land\tok{well_placed}~v~l.&&\\[3mm]
&\CoqLemma~\tok{pebbling_construction}~t:&&\\
&\quad\tok{vertex_pebbling_bound}~G~(\tok{length}~\tok{pebbles})~t\implies&&\\
&\quad\exists~(l:\coqlist~A),~l\submseteq\tok{pebbles}\land\tok{well_placed}~t~l.
\end{flalign*}
\end{tcolorbox}
\caption{A result analogous to Lemma \ref{lem:pebbling_construction}.}
\end{figure}

\pagebreak



\end{document}